\documentclass{amsart}
\usepackage{amsfonts}
\usepackage{amsmath}
\usepackage{amssymb}
\usepackage{amsthm}
\usepackage[english]{babel}
\usepackage[autostyle]{csquotes}
\usepackage{hyperref}
\usepackage{principia}
\usepackage{soul}

\MakeOuterQuote{"}

\newcommand{\keyw}[1]{\textit{Keywords---} #1}

\title{On the arithmetization of syntax}
\author{Stephen Boyce}\thanks{This paper has benefited from use of the principia.sty package.}

\newtheorem{Theorem}{Theorem}[section]
\newtheorem{Corollary}[Theorem]{Corollary}

\newtheorem{Proposition}[Theorem]{Proposition}
\newtheorem{Lemma}[Theorem]{Lemma}
\theoremstyle{definition}
\newtheorem{Definition}[Theorem]{Definition}

\newtheorem*{Lemmanon}{Lemma}
\DeclareRobustCommand{\brkbinom}{\genfrac[]{0pt}{}}

\begin{document}

\begin{abstract}
  It is generally accepted that the incompleteness of first-order number theory
  ($PA$) is established by an application of G\"odel's proof.
  This paper shows that the arithmetization of the syntax of
  $PA$ implies that the hypothesised class of $PA$ theorems is not well defined.
  This is a consequence of the fact that
  the theoremhood (or non-theoremhood)
  of a $PA$ formula is implied if the existence of
  a G\"odel number of a proof of the formula
  may be formally proved (or disproved respectively),
  using the systems own axioms and inference rules.
  Hence it is a $PA$ theorem
  that the negation ($\lnot \mathcal{G}$) of a Godel sentence $\mathcal{G}$ for $PA$
  implies  $\mathcal{G}$,
  i.e. $(\lnot \mathcal{G}) \Rightarrow \mathcal{G}$ is a $PA$
  theorem; from which, both the theoremhood and
  non-theoremhood of $\mathcal{G}$ may be established.
  The conclusion is taken as
  evidence of a failure of the devices relied upon
  for the avoidance of paradox in
  metamathematical definitions of a proof, a formal theory etc.
  The main proposition is derived using only assumptions supported by proofs already accepted in the existing literature.
\end{abstract}

\maketitle 
\begin{msc}
  03F40, 03B10, 03F30
\end{msc} \\
\keyw{Godel numberings, Classical first-order logic, Incompleteness, First-order arithmetic}

\pagenumbering{arabic}
\section{Introduction} \label{label_section_introduction}
The paper uses accepted results concerning first-order number theory
($PA$) to establish
that the of $PA$ theorems is not well defined.

The proof of this claim (Proposition \ref{propostion_theorems_not_defined})
takes only ten lines if
it is accepted that: the theoremhood (or non-theoremhood)
of a formula is implied if the existence of
a G\"odel number of a proof of the formula
may be formally proved (or disproved respectively),
using the systems own axioms and inference rules
(Proposition \ref{proposition_proof_inference}).
From this it follows that the metamathematical
definition of the class of $PA$ theorems is subject to paradox
since: it is a $PA$ theorem that the
the negation of a G\"odel sentence for
$PA$ ($\lnot \mathcal{G}$) implies $\mathcal{G}$ itself,
i.e. $(\lnot \mathcal{G}) \Rightarrow \mathcal{G}$ is a $PA$
theorem; from which, both the theoremhood and
non-theoremhood of $\mathcal{G}$ may be established.

I view the result as indicating
a failure of the devices relied upon 
for the avoidance of paradox when
the metamathematical approach is used to define a formal theory,
a formal proof etc.

Both classical logic, and the consistency of the informal
arithmetic of natural numbers are assumed. The result,
if accepted as evidence of a failure of metamathematics,
does not raise any issues with respect to either of these assumption.

With respect to the other assumptions required for the proof of
Proposition \ref{propostion_theorems_not_defined},
accepted proofs from the existing literature are
available for all of these, with one minor qualification
(discussed in \S \ref{subsection_derive_established_results}).

The reader who is familiar with this area may wish to skip directly to the proof of
Proposition \ref{propostion_theorems_not_defined}, and consider the preliminary material
only where clarification is required (concerning for example symbolic conventions
or the source for an assumed result). The only theorem which might at
first blush appear novel is in fact,
as explained in \S \ref{subsection_derive_established_results}),
supported by an accepted proof.

\section{An informal sketch}
Before entering into details,
let's briefly consider the accepted account of
the incompleteness of first-order number theory from an informal,
on some points imprecise, perspective
to assist the reader who is unfamiliar with this area.\footnote{The exposition throughout this section employs
a number of vague statements, at the cost of accuracy and correctness, to maintain brevity and accessibility.}

First-order number theory $PA$, on the metamathematical account, is the set of all sentences,
in the language of first-order arithmetic,
that are either logical truths or proper axioms corresponding to Peano's axioms or the consequences thereof.

The language of first-order arithmetic ($\mathcal{L}_A$) contains, roughly speaking:
an adequate set of propositional connectives (such as symbols for negation and
disjunction), symbols for one or more quantifiers, corresponding to the ideas of
for each and for all (that, under the intended interpretation, range over the natural numbers),
together with symbols for certain predicate and function symbols (that, under the intended interpretation,
stand for the familiar arithmetic relations and operations / functions of equality,
the successor function, addition and multiplication).

Some conventions concerning the
symbolism and syntax of $PA$ are presented in \S \ref{label_subsection_symbolism_syntax}.
For illustrative purposes, or where specifics are required,
I focus on Mendelson's \cite{mendelson2015} system $S$.
As I aim to avoid using results which are valid for only one particular version of
$PA$, the discussion of primitive symbols, axioms and inference rules
is fairly terse. These details do not alter the assumed properties of $PA$,
provided that a candidate system exhibits certain essential metamathematical properties
(\S \ref{label_section_existing_results}).

Two key properties of note that apply to all classical, metamathematical versions of $PA$ are,
somewhat confusingly, referred to as completeness and incompleteness.

By a standard result due to G\"odel \cite{godel1929}, a first-order theory such as $PA$ is complete
(if correctly formulated) in the sense that every logically valid sentence
of the language of the system ($\mathcal{L}_A$) is formally provable
using only the logical axioms and inference rules of the system.

Of course by an even more famous demonstration, also originating with G\"odel \cite{godel1931},
any (recursively defined) formal system that contains as much arithmetic as $PA$ is either (formally) inconsistent or
incomplete in the sense that it contains sentences which are neither (formally) provable nor refutable
- provided the notion of proof is restricted to what can be formally demonstrated using only the systems own axioms and inference rules.

The now standard demonstration starts with an arithmetization
of the syntax of the system.\footnote{The assertion that the symbols,
formulas, axiom(s) / inference rule, and hence proofs, of the system are recursively defined
must be true for the system
to be in scope for G\"odel's incompleteness proof.}
That is, elements of the formal system - symbols, sequences of symbols, and sequences of sequences of symbols -
are mapped into the natural numbers
in such a way that the key metamathematical properties and relations of these, such as being a formula,
a sentence or a proof, are associated with purely arithmetic properties and
relations.\footnote{More precisely, this method of G\"odel encoding a formal system, arithmetization,
was the approach used in G\"odel's original proof, though alternative methods have subsequently
been used. Shankar encodes $Z_2$ expressions in (a Lisp representation of) "$Z_2$ sets", i.e. variable-free
$Z_2$ terms (\cite{shankar1994}: 99). \'Swierczkowski similarly encodes $HF$ / $Z_2$
expressions in such $HF$ sets (\cite{swierczkowski2003}: 5).}
In this way, formulae of our system of interest, $PA$,
that are associated with arithmetic sentences or statements, are also,
indirectly, associated with statements of metamathematical facts
concerning elements of our system.
Additionally, since some sentences $\mathcal{F}$ of the system,
using numerals or other terms, predicate properties or relations of the numbers thus assigned to these
sentences themselves, a syntactically defined property of our system may be compared to a semantic phenomenon of indirect self-reference.

Since moreover the relation of being a proof may be thus arithmetized,
this implies the existence a sentence $\mathcal{G}$ which, in a sense,
diagonalizes the systems own provability relation, so that:
the assumption that the system is consistent - that contradictions, broadly speaking,
can't be proved within the system - implies that $\mathcal{G}$ is formally undecidable
(i.e. neither formally provable
nor refutable).\footnote{G\"odel's original proof shows that the assumption that the system is $\omega$-consistent
implies that $\mathcal{G}$ is formally undecidable. Following an innovation due to Rosser,
one may show that the assumption of formal consistency implies that a sentence similar to $\mathcal{G}$ is
formally undecidable (\cite{mendelson2015}: 210).}

In the decades following G\"odel's original work,
perspectives on the existence of formally undecidable propositions
have developed which reduce the issue to the question of whether certain
sets are recursive, or, alternatively, the question of
what can be determined by a Turing machine. Such material is beyond the
scope of the following discussion, which focuses instead on
the application of G\"odel's technique of arithmetizing syntax
to classical (metamathematical) systems of first-order arithmetic. 
The novel conclusion advanced below
is established by showing that 
the theoremhood (or non-theoremhood)
of a formula is implied if the existence of
a G\"odel number of a proof of the formula
may be formally proved (or disproved respectively)
using the systems own axioms and inference rules.
Hence the hypothesised class of theorems is not well defined,
since the theoremhood of a G\"odel sentence for
$PA$ implies syntactically that it is (formally) refutable,
using $PA$'s own axioms and inference rules, and vice versa.
The arithmetized metatheory of $PA$ is subject to paradox.
\section{The failure of the metamathematical notion of proof}

To establish that any standard proof of the incompleteness of
$PA$ may be extended as claimed, a number of definitions
and propositions taken from the existing literature are used. To keep matters brief,
for some established definitions and propositions I assume the generally accepted material and
cite Mendelson \cite{mendelson2015} for details. Wherever convenient I treat $PA$ and
Mendelson's first-order theory $S$ as interchangeable. The demonstration aims however to avoid exploiting any
features of $S$ that are unique by comparison with other standard accounts. It should therefore be applicable,
with appropriate changes, to any standard account of $PA$. To steer clear of an inessential
digression I use "recursive" throughout where, in several places,
the more restrictive notion of "primitive recursive" could be used (see  \cite{mendelson2015}: \S 3.3).

\subsection{Symbolism and syntax} \label{label_subsection_symbolism_syntax}
Most symbols involving a detailed definition are introduced in \S \ref{label_section_existing_results}, however
some of the more basic conventions concerning symbolism and syntax used below may be summarised as follows.

For the object language $PA$ I use
Mendelson's (\cite{mendelson2015}: Chapter 2-3) symbolism for the first-order theory
$S$.\footnote{For some symbols, the informal / unofficial symbols of $S$ are used in place of the
less familiar official symbols. The uses associated with the non-logical symbols are of course their intended interpretation.}
\begin{enumerate}
\item For logical constants
  $(\forall x_1)$, $(\exists x_1)$, $\lnot$, $\Rightarrow$, $\Leftrightarrow$, are used for, respectively, for all $x_1$, for some $x_1$, the negation,
  the conditional, and the biconditional.
\item A single individual constant "$0$" for zero is used.
\item A denumerable supply of indexed individual variables ($x_1, x_2, \ldots $) is used.
\item One (two-argument) predicate symbol for equality "$=$" is
  used.
\item One single argument function symbol ("$S$") for the successor function and
  two two-argument function symbols - "$+$" and "$\times$", for addition and multiplication respectively - are used.
\item Where no difficulties will arise some official brackets are generally omitted.
\end{enumerate}

Details concerning the formation and inference rules of the $PA$ or $S$ are, for the most part, omitted
since the only essential properties required for the proof are:
\begin{enumerate}
\item The formation and inference rules may be recursively defined;
\item The properties of the formal theory $S$ assumed for the purposes of the proof
  are shared by other standard, classical formulations of $PA$;
\end{enumerate}

For the metamathematical description of $S$:
\begin{enumerate}
\item Quotation is usually avoided since it should always be clear from context whether use
  or mention is intended.\footnote{Quine's quasi-quotation (\cite{quine1981}: \S 6) might have improved precision
  at certain points, were it not for the fact that some readers would
  incorrectly assume a convention for naming G\"odel numbers was in play.}
\item \label{item_syntactical_variables} Syntactical variables for
  individual variables ($x$, $y$, \ldots) and terms ($b$, $t$, \ldots ) are used.
\item $\overline{n}$ is used for the $PA$ numeral for the natural
  number $n$ (i.e. "$0$" nested within $n$ applications of the successor symbol - $0$, $S(0)$, $S(S(0))$ etc.).
\item Scripted letters, such as $\mathcal{F}$, are used as syntactical variables
  for (well-formed) formulae (wffs):
  \begin{enumerate}
  \item $\mathcal{F}(x)$ ($\mathcal{F}(x, y)$ etc.) are used for
    wffs that may contain one or more instances of $x$ (or $x$ and $y$ respectively) as a free  variable. 
  \item $\mathcal{F}(b)$ for the wff that results when all instances of the one free variable of
    $\mathcal{F}(x)$, should it exist, are replaced by the term $b$.
  \end{enumerate}
\item $PA$ proof: The metamathematical notion of a formal $PA$ proof is used.
  That is, the assertion that there exists a $PA$ proof of a wff $\mathcal{C}$
  is true iff there exists a finite sequence of wffs $\mathcal{B}_0$, $\ldots$, $\mathcal{B}_n$
  such that: $\mathcal{C}$ is $\mathcal{B}_n$ and each formula in the sequence
  is either a $PA$ axiom, or follows from one or more previous formulae in the sequence
  via a $PA$ inference rule.
\item  The metamathematical expression $\vdash_{PA} \mathcal{F}$ (or $\nvdash_{PA} \mathcal{F}$)
  asserts that $\mathcal{F}$ is (or is not respectively) a $PA$ theorem.
\item Conditional $PA$ proof: Let $\mathcal{B}_0$, $\ldots$, $\mathcal{B}_n$,  $\mathcal{C}$ all be wffs.
  $\mathcal{B}_0 \ldots \mathcal{B}_n \vdash_{PA} \mathcal{C}$ asserts that:
  if $PA$ proofs of $\mathcal{B}_0$, $\ldots$, $\mathcal{B}_n$ exist then
  there also exists a $PA$ proof of $\mathcal{C}$.
\item As the syntax of $PA$ may be arithmetized or G\"odel encoded in various ways,
  I assume throughout the following that some fixed though unspecified arithmetization is used.
  \begin{enumerate}
  \item Without loss of generality, where specifics are required I assume that Mendelson's arithmetization of
    $S$ is used for this purpose (\cite{mendelson2015}: \S 3.4).
  \item  Various symbols introduced by G\"odel in relation to the arithmetization of syntax
    will also be used, on the understanding that these are defined using the above mentioned arithmetization of $PA$.
  \end{enumerate}
\item Standard notions of (primitive) recursive number-theoretic
  functions and relations (\cite{mendelson2015}: \S 3.2-3)
  are used.
\item In addition to the metamathematical notion of a formal $PA$ proof,
  a primitive notion of proof and associated
  classical logical notions are used in metatheoretical demonstrations.
\end{enumerate}
For metamathematical reasoning more generally I follow G\"odel's \cite{godel1931} use of Hilbert's \cite{hilbertackermann1950} symbolism
for logical constants corresponding to those mentioned above:
$(x)$, $(Ex)$, $\overline{R}(x, y)$, $\rightarrow$, $\equiv$. Comparing such use with the use of
syntactical variables described above (at Item \ref{item_syntactical_variables})
it should be apparent the reader needs to consider the context of an expression to determine the meaning or range of
a metamathematical variable - in some contexts, this is restricted to the formal variables of $S$, whereas
in other contexts the same variable ranges over the natural numbers, and so on.

To expand the above comment on the formation and inference rules of $PA$, let's briefly anticipate a point
arising, concerning the content of $S$'s logical axioms, used below in the proof.
While the proof of Proposition \ref{proposition_proof_inference}
makes use of an instance of $S$'s universal instantiation ([10] A4) axiom, it is not essential, for this argument to succeed, that
some alternative version of $PA$ also includes as logical axioms all instances of formulae exhibiting the relevant form. All that is required
is that, where a formula corresponding to $(\forall x) \mathcal{F}(x, \ldots)$ is a theorem, the formula corresponding to $\mathcal{F}(b, \ldots)$ is also
a theorem, whenever $b$ is a term that is free for $x$ in the usual sense.
Hence, to reiterate, the only further details concerning $S$'s formation rules, logical / proper axioms, and inference rules
are various metamathematical properties discussed in \S \ref{label_section_existing_results}.

\subsection{Results taken from the existing literature}\label{label_section_existing_results}

The demonstration of the main result
(Proposition \ref{propostion_theorems_not_defined})
assumes two categories of established results presented
in this section.

The first category of assumed
results is made up of a number of metatheoretical propositions
about $PA$ which assert that one or more derived inference rule
applies to the system. These assert that:
whenever the conditions of application of the proposition hold,
it may be inferred that there exists a $PA$ proof of a $PA$ formula of interest.
Use of such rules herein is justified by
citing a standard proof. The proofs
generally proceed by induction on the complexity of an arbitrary $PA$ formula,
thus showing by cases that, however the formula is composed the rule in question will be satisfied.

The second category of assumed
result is composed of various definitions and propositions
concerned with the syntax of $PA$ broadly speaking,
including:
\begin{enumerate}
\item The existence of various recursive arithmetic
  functions and relations that correspond to syntactic notions
  concerning $PA$;
\item The existence of $PA$ formulae that (syntactically)
  define these functions and relations
  within $PA$.
\end{enumerate}

\subsubsection{Metatheorems concerned with $PA$ Theoremhood}
The first derived inference rule for discussion is the Tautology Theorem:\footnote{Hereon
the calligraphic script variables ($\mathcal{B}$, $\mathcal{C}$ etc.)
are restricted, unless otherwise indicated,
to wffs of $S$ / $PA$; for emphasis this will sometimes be explicitly stated.}
\begin{Proposition}\label{proposition_tautology}
  Tautology Theorem. If the (propositional) form of a $PA$ formula $\mathcal{C}$ corresponds to a tautology
  then there exists a $PA$ proof $\mathcal{B}_0$, $\ldots$, $\mathcal{B}_n$
  of $\mathcal{C}$
  in which only the propositional logical axioms and inference rules are used to
  justify the inclusion of a formula $\mathcal{B}_i$, $0 \leq i \leq n$, in the proof sequence;
  i.e. neither the proper axioms of $PA$,
  nor the quantification specific logical axioms /
  inference rules are used  to
  justify the inclusion of a formula in the sequence.
  (\cite{mendelson2015} Proposition 2.1).\footnote{For Mendelson's $S$,
  the $PA$ axioms concerned with quantification are: an axiom schema
  for universal instantiation ($A4$), an axiom schema concerned
  with the distribution of a (universal) quantifier
  through a conditional ($A5$), and a Generalization inference rule (\cite{mendelson2015} \S 2.3.1, \S2.3.3).}
\end{Proposition}
The intuitive notions of propositional form and correspondence
used Proposition \ref{proposition_tautology} are left undefined,
since detailed discussion of the syntax of the system, that is
not constant between different formulations of $PA$, is required.
The intuitive notion of correspondence used in Proposition \ref{proposition_tautology}
can be made precise using, for example, a truth-table algorithm.

The second derived inference rule to be considered is known as the Deduction Theorem.
This is a rule for establishing (conditional) $PA$ theoremhood.

For the Propositional Calculus the Theorem states that:
where $\mathcal{B}_0$, $\ldots$, $\mathcal{B}_n$, $\mathcal{C}$ are all wffs,
if ($\alpha$) the hypothesis that there exist $PA$ proofs of $\mathcal{B}_0$, $\ldots$, $\mathcal{B}_n$
implies that a $PA$ proof of $\mathcal{C}$ exists, then
($\beta$) $\mathcal{B}_0 \ldots \mathcal{B}_n \vdash_{PA} \mathcal{B}_n \Rightarrow \mathcal{C}$ holds
(\cite{mendelson2015}: Proposition 1.9).

If we assume that Generalisation is an inference rule of $PA$
- i.e.  $(\forall x_i) \mathcal{F}$ may be inferred from $\mathcal{F}$ (\cite{mendelson2015}: 64) -
then a statement of the Deduction Theorem
that is appropriate for $PA$ is as follows:\footnote{Mendelson also proves a version of
the deduction theorem for $PA$ (\cite{mendelson2015} Proposition 2.5)
that has slightly broader scope than the version stated above,
however the narrower statement is more succinct and adequate for our purpose.}
\begin{Proposition}\label{proposition_deduction}
  Deduction Theorem. Let $\tau$ be a (possibly non-empty) set of a finite number of $PA$ formulae.
  Let $\mathcal{B}$ and $\mathcal{C}$ be possibly distinct $PA$ formulae.
  If $\tau, \mathcal{B} \vdash_{PA} \mathcal{C}$ holds 
  (i.e. the hypothesis that there exists a $PA$ proof of $\mathcal{B}$ and each member of $\tau$
  implies that $\mathcal{C}$ is a $PA$ theorem), 
  then $\tau \vdash_{PA} \mathcal{B} \Rightarrow \mathcal{C}$ also holds,
  provided that there is no application of Generalisation in this deduction in which a free
  variable of $\mathcal{B}$ is quantified.(\cite{mendelson2015} Proposition 2.6).
\end{Proposition}
If $\tau$ is the empty set, then such an application of the deduction theorem establishes
$\vdash_{PA} \mathcal{B} \Rightarrow \mathcal{C}$ unconditionally.
That is, $\vdash_{PA} \mathcal{B} \Rightarrow \mathcal{C}$ follows
from $\mathcal{B} \vdash_{PA} \mathcal{C}$.
Informally, this corresponds to the idea that:
if we are able to prove some statement $\mathcal{E}$ on the supposition that $\mathcal{H}$ is true,
then we have proved the conditional ``If  $\mathcal{H}$ then  $\mathcal{E}$'' (\cite{mendelson2015} \S 1.9).

A third derived inference rule used in the formal demonstration
in \S \ref{propostion_theorems_not_defined} is known as
biconditional elimination (\cite{mendelson2015} \S 2.5). Mendelson
actually groups six related rules together under this heading,
all concerned with the situation in which it is hypothesised that
the proof of a biconditional exists (\cite{mendelson2015}: 75).
To illustrate, the form applied in justifying the inference applied at
step (3) in the proof of Proposition \ref{propostion_theorems_not_defined},
when $\mathcal{B}$ is $(\lnot \mathcal{G})$ and $\mathcal{C}$ is  $(\exists x_2) B(x_2, \overline{g})$, is
as follows:
\begin{Lemma} \label{lemma_biconditional_elimination} Biconditional elimination (\cite{mendelson2015}: 75).
  \begin{equation*}
    \mathcal{B} \Leftrightarrow \mathcal{C}, \mathcal{B} \vdash_{PA} \mathcal{C}
  \end{equation*}
\end{Lemma}
\begin{proof}
  In $S$, $\mathcal{B} \Leftrightarrow \mathcal{D}$ is an abbreviation for
  $(\mathcal{B} \Rightarrow \mathcal{D}) \land (\mathcal{D} \Rightarrow \mathcal{B})$, and
  $\mathcal{B} \land \mathcal{D}$ is in turn an abbreviation for
  $\lnot (\mathcal{B} \Rightarrow \lnot \mathcal{D})$.
  By the Tautology Theorem,
  Conjunction elimination holds for $PA$ so that, for any
  $PA$ wffs $\mathcal{B}$ and $\mathcal{C}$:
  \begin{Lemmanon}[Conjunction elimination (\cite{mendelson2015}: 33, $1.48(d)$]
    \begin{equation*}
      \vdash_{PA} (\mathcal{B} \land \mathcal{C}) \Rightarrow \mathcal{B}
    \end{equation*}
  \end{Lemmanon}
  Whilst we could also use the Tautology Theorem to obtain the Biconditional elimination
  lemma, an alternative approach to proving the Lemma, for the system $S$ using Conjunction elimination, is as follows:
  \begin{quote}
    \begin{flalign*}
      \quad (1)\  \mathcal{B} \Leftrightarrow \mathcal{C}  && \text{Hypothesis}
    \end{flalign*}
    \begin{flalign*}
      \quad (2)\ \mathcal{B} && \text{Hypothesis}
    \end{flalign*}
    \begin{flalign*}
      \quad (3)\  \vdash_{PA} \mathcal{B} \Rightarrow \mathcal{C} && \text{(1), conjunction elimination, modus ponens}
    \end{flalign*}
    \begin{flalign*}
      \quad (4)\ \vdash_{PA} \mathcal{C} && \text{(2), (3), modus ponens}
    \end{flalign*}
  \end{quote}
\end{proof}

Some observations concerning the metamathematical notion of "logical truth"
will facilitate a more efficient discussion of derived inference rules, and how they
may be justified for a system like $PA$.
\paragraph{The Metamathematical Notion of Logical Truth}\label{paragraph_metamathematical_notion_of_logical_truth}
For derived inference rules concerned with logical truths,
a justification that applies to any
standard, classical version of $PA$ can be obtained
by applying either the tautology theorem or the completeness of the calculus.

According to the metamathematical approach, the relevant notion of
a $PA$ "logical truth", or a logically valid $PA$ formula,
may be made precise in either one of two ways:
\begin{Definition}[$PA$ Logical Truths]\label{definition_pa_logical_truths}
  \begin{enumerate}
  \item A formula $\mathcal{B}$ of $PA$ is logically valid,
    and a semantic $PA$ Logical Truth, iff
    if it is true in every interpretation of $PA$ (\cite{mendelson2015}: 62).
  \item \label{item_pa_logical_truths_syntactic} A formula $\mathcal{B}$ of $PA$ is also logically valid,
    and a Syntactic $PA$ Logical Truth iff there exists a finite sequence of $PA$ wffs
    $\mathcal{A}_0$, $\mathcal{A}_1$, \ldots $\mathcal{A}_n$
    such that: (i) Each formula in the sequence
    is either a $PA$ logical axiom or follows from one or more previous formulae in the sequence
    via a $PA$ inference rule; and (ii) $\mathcal{B}$ is the last member of the sequence $\mathcal{A}_n$.
  \end{enumerate}
\end{Definition}
By the completeness theorem,
which is assumed by Definition \ref{definition_pa_logical_truths},
one and the same the same class of formulae is
determined to be logically valid / a logical truth in the language of arithmetic,
irrespective of which part of the definition is used.
\begin{Proposition}\label{proposition_completeness}
  Completeness Theorem. A $PA$ wff $\mathcal{B}$ is logically valid,
  in the sense of being true in every interpretation, if and only
  if it is a  Syntactic $PA$ Logical Truth
  (Definition \ref{definition_pa_logical_truths} \ref{item_pa_logical_truths_syntactic}),
  i.e. established as a $PA$ theorem on the basis of $PA$'s
  logical axioms and inference rules (\cite{mendelson2015} Corollary 2.19).
\end{Proposition}
With the Tautology Theorem and Completeness Theorem in hand,
we now have two general derived inference rules that
apply to any standard version of $PA$,
irrespective of variations with respect to syntax. Hence
I will now simply list four additional derived inference rules
used below which may be thus justified, before considering
some additional inferences rules that warrant
some further clarification.
\begin{Lemma} Conjunction introduction (\cite{mendelson2015}: 75).
  \begin{equation*}
    \mathcal{B}, \mathcal{C} \vdash_{PA} \mathcal{B} \land \mathcal{C}
  \end{equation*}
\end{Lemma}

\begin{Lemma} Biconditional negation (\cite{mendelson2015}: 75).
  \begin{equation*}
    \mathcal{A} \Leftrightarrow \mathcal{B} \vdash_{PA} (\lnot \mathcal{A}) \Leftrightarrow (\lnot \mathcal{B})
  \end{equation*}
\end{Lemma}

\begin{Lemma} Substitutivity of equality (\cite{mendelson2015}: 93).
  Since $S$ is a theory with equality (\cite{mendelson2015}: Corollary 3.3),
  \begin{quote}
    \ldots the following is a theorem of [$S$]
    \begin{flalign*}
      x = y \Rightarrow (\mathcal{B}(x, x) \Rightarrow \mathcal{B}(x, y)) \tag{A7} &&
    \end{flalign*}
    where $x$ and $y$ are any variables,
    $\mathcal{B}(x, x)$ is any wf, and $\mathcal{B}(x, y)$
    arises from $\mathcal{B}(x, x)$ by replacing some, but not necessarily all, free occurrences of $x$
    by $y$, with the proviso that $y$ is free for $x$ in $\mathcal{B}(x, x)$ \ldots
    (\cite{mendelson2015}: 93, modified at square brackets, "$K$" changed to "$S$").
  \end{quote}
\end{Lemma}

\begin{Lemmanon} Equivalence theorem (\cite{mendelson2015}: Proposition 2.9(a)).
  If the following two conditions hold:
  \begin{enumerate}
  \item $\mathcal{C}$ and $\mathcal{D}$ are alike excepting that at certain locations: \newline
    $\mathcal{C}$ contains an occurrence $\mathcal{A}$ \newline
    where $\mathcal{D}$ contains an occurrence $\mathcal{B}$;
  \item $y_1$, \ldots, $y_k$ comprises the entire list of variables $y$ such that:\newline
    $y$ is a free variables of $\mathcal{A}$ that is bound in some occurrence of $\mathcal{A}$ in $\mathcal{C}$,\newline
    or $y$ is a free variables of $\mathcal{B}$ that is bound in some occurrence of $\mathcal{B}$ in $\mathcal{D}$;
  \end{enumerate}
  then the following theorem also obtains:
  \begin{equation*}
    \vdash_{PA} [(\forall y_1) \ldots (\forall y_k) (\mathcal{A} \Leftrightarrow \mathcal{B})] \Rightarrow (\mathcal{C} \Leftrightarrow \mathcal{D})
  \end{equation*}
\end{Lemmanon}

A number of other derived inference rules used below require
a small amount of discussion. The first of these is
described below as "transitivity of conditional proof" (e.g. step (5) in the proof of
Proposition \ref{propostion_theorems_not_defined}).
\paragraph{Transitivity of Conditional Proof}\label{paragraph_transitivity_of_conditional_proof}
For the system $S$, Mendelson
does not give this inference principle a name, but introduces it as a property of the
notion of consequence. The principle is stated in a slightly different form
to that used below and justified as follows (where, for present purposes,
$\Delta$ and $\tau$ may be taken to
be finite sets of $S$ wffs and $\vdash$ may be taken as $\vdash_{PA}$):
\begin{quote}
  3. If $\Delta \vdash \mathcal{C}$, and for each $\mathcal{B}$ in $\Delta$,
  $\tau \vdash \mathcal{B}$ then $\tau \vdash \mathcal{C}$. \newline
  \ldots if $\mathcal{C}$ is provable from premisses in $\Delta$,
  and each premiss in $\Delta$ is provable from premisses in
  $\tau$, then ${C}$ is provable from premisses in $\tau$. (\cite{mendelson2015}: 28).
\end{quote}
Use of the rule, stated in this way, is generally indicated by lines of the following form
(e.g. proof of Lemma 1.11(e) line 8, \cite{mendelson2015}: 33):
\begin{flalign*}
  \quad n.\  \quad \quad \mathcal{B}  \vdash \mathcal{C}\ && 1 - (n{-}1)
\end{flalign*}
Mendelson's statement of the rule
and intuitive justification are perfectly fine.
If the rule is stated in a slightly different form
however a more formal justification suggests itself.
\begin{Lemmanon}[Transitivity of Conditional Proof]
  Let $\tau$ be a (possibly non-empty) set of a finite number of $PA$ formulae.
  Let $\mathcal{B}_0$, \ldots $\mathcal{B}_n$ and $\mathcal{C}$ be possibly distinct $PA$ formulae. Assume
  that each of the following assertions of conditional theoremhood, $(0)-(n)$, are true:
  \begin{flalign*}
    \quad (0)\ \tau,\ \mathcal{B}_0 \vdash_{PA} \mathcal{B}_1
  \end{flalign*}
  \begin{flalign*}
    \quad (1)\  \tau,\ \mathcal{B}_1 \vdash_{PA} \mathcal{B}_2
  \end{flalign*}
  \begin{flalign*}
    \quad \ldots
  \end{flalign*}
  \begin{flalign*}
    \quad (n-1)\ \tau,\ \mathcal{B}_{n-1} \vdash_{PA} \mathcal{B}_n
  \end{flalign*}
  \begin{flalign*}
    \quad (n)\ \tau,\ \mathcal{B}_n \vdash_{PA} \mathcal{C}
  \end{flalign*}
  Then the following assertion of conditional theoremhood is also true:
  \begin{equation*}
    \quad (n+1)\ \tau,\ \mathcal{B}_0 \vdash_{PA} \mathcal{C}
  \end{equation*}
\end{Lemmanon}
\begin{proof}
  If one or more of the hypotheses,
  $\tau,\ \mathcal{B}_i \vdash_{PA} \mathcal{B}_{i+1}$,
  or   $\tau,\ \mathcal{B}_n \vdash_{PA} \mathcal{C}$,
  are false it follows
  immediately that the conditional is true.
  To see if the Lemma can fail, let's assume that
  each of the hypotheses $(0)-(n)$ are true.
  Let each hypothesis then be
  associated with a finite sequence of $S$ wffs
  ($\mathcal{A}_{i_0}$, $\mathcal{A}_{i_1}$, \ldots $\mathcal{A}_{i_j}$)
  which corresponds to the
  hypothesised conditional proof of $\mathcal{B}_{i+1}$,
  or $\mathcal{C}$, premised on the hypothesised theoremhood
  of $\mathcal{B}_i$ and each member of $\tau$,
  so that:
  \begin{enumerate}
  \item For each formula in the union of $\tau$, \{$\mathcal{B}_i$\},
    and \{$\mathcal{B}_{i+1}$\} or  \{$\mathcal{C}$\},
    there exists an $m$, for $0 \leq m \leq i_j$,
    such that this formula is $\mathcal{A}_{i_m}$ in the sequence.
  \item $\mathcal{B}_{i+1}$ or $\mathcal{C}$ is the last member of the sequence $\mathcal{A}_{i_j}$,
    although $\mathcal{B}_{i+1}$ and $\mathcal{C}$, as
    well as $\mathcal{B}_i$ and the members of $\tau$, may occur more than once in the sequence.
  \item The sequence corresponds to the hypothesised proof of $\mathcal{B}_{i+1}$
    or $\mathcal{C}$, conditional upon
    the hypothesised (conjoint) theoremhood of $\mathcal{B}_i$ and each member of $\tau$,
    so that each formula in the sequence
    is either ($\alpha$) equal to $\mathcal{B}_i$; or ($\beta$) a member of $\tau$; or
    ($\gamma$) a $PA$ axiom; or ($\delta$) follows from one or more previous formulae in the sequence
    via a $PA$ inference rule.
  \end{enumerate}
  The third of these items is less restrictive than the corresponding
  requirement for an (unconditional metamathematical) proof - via the clauses $\alpha$ and $\beta$ -
  so that an assertion of conditional theoremhood, such as ($\tau,\ \mathcal{B}_0 \vdash_{PA} \mathcal{C}$),
  may, if appropriate, be determined to be true
  even if the hypothesised $PA$ proofs of $\mathcal{B}_0$ and of each member of $\tau$ do not exist.
  For example, if we assume that $PA$ is consistent, then for no wff $\mathcal{B}$
  is there a $PA$ proof of the formula $\mathcal{B} \land (\lnot \mathcal{B})$;
  yet, where $\mathcal{B}$ and $\mathcal{C}$ are arbitrary wffs,
  the assertion of conditional theoremhood ($\mathcal{B} \land (\lnot \mathcal{B}) \vdash_{PA} \mathcal{C}$) is true.
  
  If now we form the sequence of sequences of $S$ wffs corresponding to this hypothesised sequence of
  conditional proofs,
  then the existence of the required conditional $S$ proof, $\tau,\ \mathcal{B}_0 \vdash_{PA} \mathcal{C}$,
  is established (where
  the expression to the left of the turnstile denotes the set of formulae corresponding to the sequence just mentioned):
  \begin{flalign*}
    \begin{split}
      \quad \mathcal{A}_{0_0}, \mathcal{A}_{0_1},\ \ldots,\ \mathcal{A}_{0_j}\ \mathcal{A}_{1_0}, \mathcal{A}_{1_1},\ \ldots,\ \mathcal{A}_{1_k},\ \ldots,\ \\
      \mathcal{A}_{(n+1)_0}, \mathcal{A}_{(n+1)_1},\ \ldots,\ \mathcal{A}_{(n+1)_l}\ \vdash \mathcal{C} &
    \end{split}
  \end{flalign*}
  Hence the consequence is true and failure of the Lemma is excluded.
\end{proof}
A second pair of derived inference rules used below concerns the
introduction and elimination of defined expressions.
\paragraph{Introduction and elimination of definitions}
The system $S$, like many modern systems of logic, makes use of
definitions. Four used throughout the following discussion are as follows.
\begin{Definition}\label{definition_abbreviation}
  Where $\mathcal{B}$ and $\mathcal{C}$ are any wffs and $x$ is any variable,
  the expression on the left of "$=_{\text{def}}$" shall be used as an abbreviation
  for the expression appearing on the right (adapted from \cite{mendelson2015}:29, 43, 98, \cite{pm1910v1}):
  \begin{quote}
    \vspace{-4ex}
    \begin{flalign*}
      \quad (\mathcal{B} \land \mathcal{C}) =_{\text{def}} \lnot (\mathcal{B} \Rightarrow \lnot \mathcal{C}) \tag{D1} &&
    \end{flalign*}
    \begin{flalign*}
      \quad (\mathcal{B} \lor \mathcal{C}) =_{\text{def}} (\lnot \mathcal{B}) \Rightarrow \mathcal{C} \tag{D2} &&
    \end{flalign*}
    \begin{flalign*}
      \quad (\mathcal{B} \Leftrightarrow \mathcal{C}) =_{\text{def}} (\mathcal{B} \Rightarrow \mathcal{C}) \land (\mathcal{C} \Rightarrow \mathcal{B}) \tag{D3} &&
    \end{flalign*}
    \begin{flalign*}
      \quad ((\exists x) \mathcal{B})  =_{\text{def}} \{\lnot [(\forall x) (\lnot \mathcal{B})]\} \tag{D4} &&
    \end{flalign*}
    \begin{flalign*}
      ((\exists_1 x) \mathcal{B})  =_{\text{def}} \{(\exists x) \mathcal{B}(x)\ & \land  \quad  \quad  \tag{D5} \\
      & (\forall x)(\forall y) [\mathcal{B}(x) \land \mathcal{B}(y) \Rightarrow x = y]\}  \quad \quad \quad \quad 
    \end{flalign*}
  \end{quote}
\end{Definition}
The reader will observe that some conventions regarding the introduction and elimination or parentheses
are assumed at this point. The subtleties of the introduction and elimination of defined expressions,
beyond the observation that the defined expressions are to be treated as abbreviations (\cite{mendelson2015}:29),
are apparently assumed knowledge. The following convention, borrowed from \emph{Principia} \cite{pm1910v1},
recommends itself - the description being quite informal, though clearly amenable to a more precise statement
if required.

The official $S$ symbols, formulae, proofs etc. do not include any defined symbols. Thus, to obtain,
for example, an official wf from an unofficial expression in which one or more defined symbols appear
each such defined symbol must be eliminated. The order for the elimination of defined terms
shall be the reverse of the typographic order ($D1$, $D2$, \ldots.) in which the abbreviations are defined.
For example, if abbreviations of all kinds occur in the the expression $\mathcal{B}$, abbreviations introduced
using $D5$ should be eliminated prior to the elimination of abbreviations introduced
using $D4$ etc.

When eliminating instances of the $nth$ such defined term (e.g. the fourth term being $D4$), one commences
with the instance of least scope, and, using the equivalence as provided by the above definitions,
replaces this instance involving a defined term with the equivalent expression that includes one
less instance of a defined symbol / sequence of symbols. This process is repeated, continuing with the
next remaining abbreviation of this kind that has the least scope, until all instances of this abbreviation are thus removed.
The process is then repeated with each remaining kind of abbreviation, until an official expression involving
no abbreviations is obtained.

Since a precise description of the process of replacing abbreviations is quite complicated,
and is not corroborated by an official account, an alternative process is followed below. Firstly,
the following lemma is used.
\begin{Lemma}[Equivalence by Definition] \label{lemma_abbreviation}
  Where $\mathcal{B}$ and $\mathcal{C}$:
  \begin{enumerate}
  \item are like expressions built from official $S$ symbols
    except for involving one or possibly more $S$ abbreviations
    using the definitions at Definition \ref{definition_abbreviation}; 
  \item are both transformed to $S$ wffs, once all abbreviations are eliminated, using
    the definitions at Definition \ref{definition_abbreviation};
  \item are both transformed to the same, possibly distinct,
    $S$ formula, $\mathcal{D}$,
    once all abbreviations are eliminated using the definitions at Definition \ref{definition_abbreviation};
  \end{enumerate}
  Then the following theorem holds:
  \begin{equation*}
    \vdash_{PA} \mathcal{B} \Leftrightarrow \mathcal{C}  \tag{\ref{lemma_abbreviation}.1}
  \end{equation*}
  In Equation \ref{lemma_abbreviation}.1, and wherever convenient hereon,
  "Metamathematical statements about terms and formulas of the system
  are \ldots to be understood to refer to the un-abbreviated expressions" (\cite{kleene1962}: 75).
  Hence, the metamathematical variables in Equation \ref{lemma_abbreviation}.1
  are taken to be names for the wffs obtained when all abbreviations are
  eliminated.
\end{Lemma}
\begin{proof}
  From the truth of all three hypotheses we have that both $\mathcal{B}$ and $\mathcal{C}$ are
  transformed (applying the definitions at Definition \ref{definition_abbreviation} to eliminate
  zero or more abbreviations) to a common, possibly distinct,
  $S$ formula, $\mathcal{D}$. (If one or more of these hypotheses is false then the conditional is
  trivially true and we are done.)
  For the proof of the theorem we need only establish
  that: where $\mathcal{D}$ is any $S$ wff, $\mathcal{D} \Leftrightarrow \mathcal{D}$
  is an $S$ theorem. This may be established as follows.
  Let's note firstly that, by the Tautology Theorem,
  the law of Identity holds for $PA$ so that, for any
  $PA$ formula $\mathcal{B}$:
  \begin{Lemma} \label{lemma_identity} Identity (\cite{mendelson2015}: Lemma 1.8).
    \begin{equation*}
      \vdash_{PA} \mathcal{B} \Rightarrow \mathcal{B}
    \end{equation*}
  \end{Lemma}
  Hence Equation \ref{lemma_abbreviation}.1 may be proved as follows:
  \begin{quote}
    \begin{multline*}
      \quad (1)\ \vdash_{PA} \mathcal{D} \Rightarrow \mathcal{D}\  \\
      \text{identity}
    \end{multline*}
    \begin{multline*}
      \quad (2)\ (\mathcal{D} \Rightarrow \mathcal{D}), (\mathcal{D} \Rightarrow \mathcal{D})\  
      \vdash_{PA} (\mathcal{D} \Rightarrow \mathcal{D}) \land (\mathcal{D} \Rightarrow \mathcal{D})\  \\
      \text{(1), conjunction introduction}
    \end{multline*}
  \end{quote}
\end{proof}
The proof of Lemma \ref{lemma_abbreviation}, hereon referred to as "Equivalence by Definition"
where convenient, seems to pull something out of thin air.
The substantive metamathematical propositions assumed -
that all such expressions ($\mathcal{B}$ and $\mathcal{C}$)
can be so transformed into $S$ wffs,
that each can be so transformed into a unique $S$ wffs,
etc. - are all smuggled into the hypotheses
of the Lemma without being independently established.

For an alternative approach that brings
the proof of such matters into the foreground,
we might use the concept of an
"extension by definitions" as follows.

\paragraph{Extension by definitions}\label{paragraph_extension_by_definitions}
Let $S'$ be the new formal theory 
that that is the same as $S$ except that:
\begin{enumerate}
\item The symbols defined at (D1)-(D4) - 
  i.e. $\land$, $\lor$, $\Leftrightarrow$
  and $(\exists x)$ - are added as primitive symbols;
\item Each instance of (D1)-(D4) is added as
  a primitive axiom of $S'$ (defining
  the symbol for each context of use).
\end{enumerate}
We might then prove, by induction on the complexity
of an arbitrary $S'$ formula,
that $S'$ is an "extension by definitions" of $S$
in the following sense.
\begin{Definition}[Extension by Definitions]\label{definition_extension_by_definitions}
  $S'$ is an "extension by definitions" of $S$ iff:
  For each formula $\mathcal{A}$ of $S'$
  involving one or more instances of the newly
  added symbols ($\land$, $\lor$, $\Leftrightarrow$
  and $(\exists x)$), there exists a corresponding
  formula  $\mathcal{A}^{*}$ of $S$ not involving these symbols
  (the translation of  $\mathcal{A}$ into $S$)
  such that:
  \begin{equation*}
    \vdash_{S} \mathcal{A}^{*} \leftrightarrow\ \vdash_{S'} \mathcal{A}\ \text{(\cite{shoenfield1967}: 60, paraphrased)}
  \end{equation*}
\end{Definition}
Since however the latter approach would require
more detailed discussion of the syntax of $PA$,
and a proliferation, as we progress, of propositions concerning the
translation of statements between the two formal theories,
I will hereon simply use the Equivalence by Definition Lemma.
This leads back to the question of the status of the Equivalence by Definition Lemma.

\paragraph{Proof of Equivalence by Definition}\label{paragraph_proof_of_equivalence_by_definition}
The reader looking for an independent proof of
the Equivalence by Definition Lemma in the literature may consider,
in addition to Kleene's remarks cited above,
Church's comments (\cite{church1956}: \S 11),
or Quine's characteristically pithy aphorism:
\begin{quote}
  To define a sign is to show how to avoid it. (\cite{quine1981}: 47)
\end{quote}

As Russell might have observed, it's proof by metamathematical
definition as far as the eye can see. Since the Equivalence by Definition Lemma
is both the accepted result, as defended in the above-cited cursory remarks,
and is technically correct it is used hereon to deal with abbreviations
as needed. It can be seen however that the accepted approach
essentially kicks the problem upstairs. That
is, for example, for arbitrary wffs  $\mathcal{B}$ and $\mathcal{C}$,
it is correct in some contexts to view both $(\mathcal{B} \land \mathcal{C})$
and $\lnot (\mathcal{B} \Rightarrow \lnot \mathcal{C})$ as different
metamathematical names for the same object-language formula; in the
latter case, the object-language expression itself being used in the
metalanguage as a name for instances of this very same
expression occurring in the object language. The substantive
problems still remain - of showing that the appropriate rules are being
applied to eliminate abbreviations from metamathematical names for wff etc.
The Equivalence by Definition Lemma in effect places the burden upon the user
to correctly identify and apply the transformation rules required
to eliminate abbreviations from metamathematical names for object-language formulae.

A final derived inference rules used below
that is not discussed in Mendelson (\cite{mendelson2015})
is referred to as Conditional Identity.
\paragraph{Identity}
The second line of the proof of
Proposition \ref{propostion_theorems_not_defined}
is justified by the inference rule "Conditional Identity":
\begin{Lemmanon}[Conditional Identity]
  \begin{equation}\label{equation_conditional_identity}
    \mathcal{B} \vdash_{PA} \mathcal{B}
  \end{equation}
\end{Lemmanon}
This rule is introduced purely to emphasise that at this point the proof is conditional on the
hypothesis that an $S$ proof of $\mathcal{B}$ exists.
It is of course completely redundant, since this in fact is the Hypothesis at (1) in the proof of
Proposition \ref{propostion_theorems_not_defined} and the meaning of writing
$\mathcal{B}$ to the left of the turnstile in \ref{equation_conditional_identity};
however it can be formally justified by the Identity lemma
mentioned above ($\vdash_{PA} \mathcal{B} \Rightarrow \mathcal{B}$, (\cite{mendelson2015}: Lemma 1.8).

We can now consider the assumed
results concerned more directly with the syntax of $PA$.

\subsubsection{Propositions concerned with the syntax of $PA$}
For several of the following propositions it is helpful to introduce an alternative symbolism for the operation of
substituting terms for free variables in $PA$ formulae
which shows both the variable and term explicitly:
\begin{Definition}\label{definition_informal_substitution}
  $\mathrm{Subst}\ \left[ \mathcal{B}\  \genfrac{}{}{0pt}{0}{x}{t} \right]$ (\cite{godel1931}: 155).
  Where $\mathcal{B}$ is a $PA$ formula,
  $x$ is a $PA$ variable,
  $t$ is a $PA$ term is free for $x$ in  $\mathcal{B}$,
  $\mathrm{Subst}\ \left[ \mathcal{B}\  \genfrac{}{}{0pt}{0}{x}{t} \right]$ shall be the $PA$ formula that results when
  $t$ is substituted for all free occurrences of $x$ in  $\mathcal{B}$.
  (If there are no free occurrences of $x$ in  $\mathcal{B}$ then
  $\mathrm{Subst}\ \left[ \mathcal{B}\  \genfrac{}{}{0pt}{0}{x}{t} \right]$ is just $\mathcal{B}$.)
\end{Definition}
The following four properties of $PA$ used in the proof are ultimately based on
the fact that $PA$ is a recursively axiomatised theory, and hence, the metamathematical
definitions of syntactical properties and relations of entities associated with the
formal theory $PA$ may be associated, via a G\"odel mapping, with recursive arithmetic
functions and relations. 

Firstly, the property of
being a $PA$ numeral $\overline{n}$ may be associated with a recursive function of natural numbers $Z(n)$.
\begin{Proposition}\label{proposition_arithmetized_numeral}
  $Z(n)$.
  There exists a recursive, arithmetized \textsc{numeral} function $Z(n)$ (\cite{godel1931}: 165)
  for $PA$, such that for all pairs of natural numbers $x$, $y$:\newline
  $y = Z(x)$ holds iff $y$ is the G\"odel number of the $PA$ numeral
  $\overline{x}$ for the number $x$ (cf. \cite{mendelson2015} Proposition 3.27(17)).
\end{Proposition}
Secondly, the metamathematical operation of substituting a term for a free variable in a $PA$ formula,
where the term is free for the variable in this formula,
may likewise be associated with a recursive function of natural numbers.
\begin{Proposition}\label{proposition_arithmetized_substitution}
  $\mathrm{Sb}\ \left[ f\ \genfrac{}{}{0pt}{0}{x}{t} \right]$ (\cite{godel1931}: 167).
  There exists a recursive, arithmetized substitution function for $PA$
  such that, for all quadruples of natural numbers $f$, $x_y$, $t_w$, $u$:
  where $f$ is the G\"odel number of a $PA$ formula $\mathcal{F}$,
  $x_y$ is the G\"odel number of a $PA$ variable $y$,
  $t_w$ is the G\"odel number of a $PA$ term $w$
  that is free for $y$ in $\mathcal{F}$,
  $u = \mathrm{Sb}\ \left[ f\ \genfrac{}{}{0pt}{0}{x_y}{t_w} \right]$
  holds iff $u$ is the G\"odel number of
  $\mathrm{Subst}\ \left[ \mathcal{F}\  \genfrac{}{}{0pt}{0}{y}{w} \right]$
  (cf. \cite{mendelson2015} Proposition 3.26(10)).
\end{Proposition}
Borrowing again from (\cite{godel1931}: 167),
$\mathrm{Sb}\ \left[ b\ \genfrac{}{}{0pt}{0}{x}{t} \genfrac{}{}{0pt}{0}{y}{u} \right]$
shall be used instead of \newline
$\mathrm{Sb}\ \left[ \mathrm{Sb}\ \left[ b\ \genfrac{}{}{0pt}{0}{x}{t} \right] \genfrac{}{}{0pt}{0}{y}{u} \right]$;
likewise with respect to $\mathrm{Subst}\ \left[ \mathcal{B}\ \genfrac{}{}{0pt}{0}{x}{t} \genfrac{}{}{0pt}{0}{y}{u} \right]$.

Thirdly, we may also define a recursive Diagonal Function, which,
when applied to the G\"odel number of a PA class sign
$\mathcal{B}(x_1)$ in which the first $PA$ variable is the only free variable,
returns the G\"odel number of the formula $\mathcal{B}(\overline{b})$ that results when the
numeral $\overline{b}$ for the G\"odel number of $\mathcal{B}(x_1)$ itself
is substituted for every free occurrence of $x_1$ in  $\mathcal{B}(x_1)$
\begin{Corollary}\label{corollary_diagonal_function}
  $\mathrm{Sb}\ \left[ f\ \genfrac{}{}{0pt}{0}{v_1}{Z(f)} \right]$.
  There exists a recursive, arithmetized diagonal function for $PA$
  such that for all pairs of natural numbers $f$, $u$:
  where $f$ is the G\"odel number of a $PA$ formula $\mathcal{F}$ with one free variable $x_1$,\newline
  $v_1$ is the G\"odel number of the $PA$ variable $x_1$, 
  $u = \mathrm{Sb}\ f\ {\brkbinom{v_1}{Z(f)}}$ holds iff
  $u$ is the G\"odel number of
  $\mathrm{Subst}\ \left[ \mathcal{F}\  \genfrac{}{}{0pt}{0}{x_1}{\overline{f}} \right]$
  (cf. \cite{mendelson2015} Proposition 3.27(19)).
\end{Corollary}

As the diagonal function is recursive, we may infer the following:

\begin{Lemma}[$\mathcal{D}(x_1, x_2)$ The Diagonal Relation Sign]\label{lemma_diagonal_relation_sign}
  Since the the diagonal function is recursive,
  there exists a $PA$ relation sign, $\mathcal{D}(x_1, x_2)$, that syntactically defines this function
  in $PA$ so that (\cite{mendelson2015}: 170):
  \begin{enumerate}
  \item For all pairs of natural numbers $x$, $y$:% \newline
    \begin{flalign*}
      \quad \{y\ = \mathrm{Sb}\ \left[ x\ \genfrac{}{}{0pt}{0}{v_1}{Z(x)} \right]\} \rightarrow\ \vdash_{PA} \mathcal{D}(\overline{x}, \overline{y}) && \tag{\ref{lemma_diagonal_relation_sign}.1}
    \end{flalign*}
    \begin{flalign*}
      \quad \vdash_{PA} (\exists_{1} z) \mathcal{D}(\overline{x}, z) && \tag{\ref{lemma_diagonal_relation_sign}.2}
    \end{flalign*}
  \item $\mathcal{D}(x_1, x_2)$ moreover "strongly represents" (\cite{mendelson2015}: 170) the diagonal function in $PA$,
    so that:
    \begin{flalign*}
      \quad \vdash_{PA} (\exists_{1} z) \mathcal{D}(x_1, z) && \tag{\ref{lemma_diagonal_relation_sign}.3}
    \end{flalign*}
  \end{enumerate}
\end{Lemma}

The first two parts of the Diagonal Relation Sign Lemma / Lemma \ref{lemma_diagonal_relation_sign})
follow from the fact that the diagonal function is recursive
(Corollary \ref{corollary_diagonal_function}) and
every recursive function is syntactically defined (or "representable")
in $S$ in the above sense  (\cite{mendelson2015}: Proposition 3.24).
The third part follows from the fact that every recursive function of natural numbers that is
syntactically defined in $PA$ in the sense
of the first two parts is "strongly representable"
in $PA$ (\cite{mendelson2015}: Proposition 3.12)
in the sense of the third part.

Since the metamathematical relation that holds between
a finite sequence of $PA$ formulae ($\mathcal{B}_0$, $\ldots$, $\mathcal{B}_n$)
and a possibly distinct $PA$ formula $\mathcal{C}$ when the former is a
$PA$ proof of the latter,
corresponds to a recursive arithmetic relation of natural numbers:
\begin{Proposition}\label{proposition_arithmetized_proof}
  $x\ \mathrm{B}\ y$.
  There exists a recursive, arithmetized proof relation for $PA$,
  $x\ \mathrm{B}\ y$ (\cite{godel1931}: 171),
  such that for all pairs of natural numbers $x$, $y$: 
  $x\ \mathrm{B}\ y$ holds iff $x$ is the G\"odel number of
  a proof in $PA$ of a $PA$ formula with G\"odel number $y$ (cf. \cite{mendelson2015} Proposition 3.28).
\end{Proposition}

We now consider a couple of propositions based on the facts that, firstly,
each of the metamathematical functions / relations introduced in the above propositions are recursive,
and secondly every recursive function / relation is (numeral-wise) syntactically definable,
in the required sense, in
$PA$  (cf. \cite{mendelson2015} Proposition 3.24 \cite{mendelson2015} Corollary 3.25).
Hence, both the relation $x\ \mathrm{B}\ y$ and
the diagonal function for $PA$ may be defined in $PA$ in the required proof-theoretic sense.
\begin{Proposition}\label{proposition_proof_relation_sign}
  There exists a $PA$ relation sign $B(x, y)$ that syntactically defines
  the relation of natural numbers $x\ \mathrm{B}\ y$ in $PA$; that is
  for all pairs of natural numbers $x$, $y$ (cf. \cite{mendelson2015} Proposition 3.36):
  \begin{quote}
    \vspace{-4ex}
    \begin{flalign*}
      \quad x\ \mathrm{B}\ y \rightarrow\ \vdash_{PA} B(\overline{x}, \overline{y}) \tag{\ref{proposition_proof_relation_sign}.1} &&
    \end{flalign*}
    \begin{flalign*}
      \quad \overline{x\ \mathrm{B}\ y} \rightarrow\ \vdash_{PA} \lnot B(\overline{x}, \overline{y}) \tag{\ref{proposition_proof_relation_sign}.2} &&
    \end{flalign*}
  \end{quote}
\end{Proposition}
Since the diagonal function is syntactically definable in
$PA$ (cf. \cite{mendelson2015} Proposition 3.27 (19)),
we have by an application of the Diagonalization Lemma to $PA$
(cf. \cite{mendelson2015} Proposition 3.34):
\begin{Proposition} \label{proposition_godel_sentence}
  There exists a $PA$-sentence $\mathcal{G}$, with G\"odel number $g$ 
  such that (cf. \cite{mendelson2015}: 208, Equation $\$$):
  \begin{equation*}
    \vdash_{PA} \mathcal{G} \Leftrightarrow (\forall x_2) [\lnot B(x_2, \overline{g})] \tag{\ref{proposition_godel_sentence}.1}
  \end{equation*}
\end{Proposition}
\begin{proof}
  Due to the manner in which the Diagonalization Lemma is cited above,
  it appears, contrary to fact, as though Proposition \ref{proposition_godel_sentence}
  is non-constructive. To verify this point it may be helpful to rehearse
  the now standard construction of $\mathcal{G}$ using the material at hand. The reader
  who is happy to accept Proposition \ref{proposition_godel_sentence}
  as an established result may skip this proof without loss.\footnote{The construction presented
  below is mostly a close paraphrase of Mendelson (\cite{mendelson2015}: \S 3.5),
  with a small number of minor adjustments. The adjustments
  are made as Mendelson's approach establishes the existence of the required
  G\"odel sentence for any first-order theory $K$ (with equality)
  in the language of arithmetic in scope  for
  G\"odel's First Incompleteness Theorem.
  If any errors have crept in as a result of adjusting the demonstration to focus
  specifically on the case of $S$ / $PA$ these are naturally my own.}

  Some additional symbolism will assist at this point (paraphrased with some minor modifications from (\cite{mendelson2015}: \S 3.5)):
  \begin{enumerate}
  \item Let $z$ be the first $S$ variable not used in either
    $\mathcal{D}(x_1, x_2)$ or $B(x_2, x_1)$;
  \item Let  $\mathcal{D}'(x_1, z)$ be the formula that results
    when $z$ is substituted for all free occurrences of $x_2$ in $\mathcal{D}(x_1, x_2)$;
  \item Let  $\mathcal{B}'(x_2, z)$ be the formula that results
    when $z$ is substituted for all free occurrences of $x_1$ in $B(x_2, x_1)$;
  \item \label{item_definition_m} Let $m$ the G\"odel number of the following wff with one free variable $x_1$ (hereon $\mathcal{C}(x_1)$
    where convenient):
    \begin{equation*}
      \quad \quad (\forall z) \{\mathcal{D}'(x_1, z) \Rightarrow (\forall x_2) [\lnot B'(x_2, z)]\} \tag{\ref{proposition_godel_sentence}.3}
    \end{equation*}
  \item \label{item_definition_g} Let $g$ be the G\"odel number of $\mathcal{C}(\overline{m})$, so that:
    \begin{equation*}
      \quad \quad g\ = \mathrm{Sb}\ \left[ m\ \genfrac{}{}{0pt}{0}{v_1}{Z(m)} \right] \tag{\ref{proposition_godel_sentence}.4}
    \end{equation*}
  \item \label{item_definition_wff_g} Let $\mathcal{G}$ be the above mentioned (Item \ref{item_definition_g})
    $PA$ sentence with G\"odel number $g$.
    $\mathcal{G}$ thus defined may also be viewed as a metamathematical abbreviation, in the previously explained sense
    (Equivalence by definition), for the $PA$ sentence that may also be referred to by two other metamathematical abbreviations
    introduced just above. Firstly, from Equation \ref{proposition_godel_sentence}.4 we have:
    \begin{equation*}
      \quad \quad \mathcal{G} =_{df} \mathrm{Subst}\ \left[ \mathcal{C}\  \genfrac{}{}{0pt}{0}{x_1}{\overline{m}} \right] \tag{\ref{proposition_godel_sentence}.5}
    \end{equation*}
    Secondly, substituting the alternative abbreviation for the $PA$ wff named by $\mathcal{C}(x_1)$ provided at Equation \ref{proposition_godel_sentence}.3
    we have:
    \begin{equation*}
      \quad \quad \mathcal{G} =_{df} (\forall z) \{\mathcal{D}'(\overline{m}, z) \Rightarrow (\forall x_2) [\lnot B'(x_2, z)]\} \tag{\ref{proposition_godel_sentence}.6}
    \end{equation*}
  \end{enumerate}
  That $\mathcal{G}$ thus defined is a $PA$ sentence satisfying (\ref{proposition_godel_sentence}.1) may be
  established as follows. Note firstly that, through simple substitution of the definitions just provided,
  the result to be proved may be restated as follows:
  \begin{equation*}
    \vdash_{PA} (\forall z) \{\mathcal{D}'(\overline{m}, z) \Rightarrow (\forall x_2) [\lnot B'(x_2, z)]\} \Leftrightarrow (\forall x_2) [\lnot B(x_2, \overline{g})] \tag{\ref{proposition_godel_sentence}.7}
  \end{equation*}
  From the definition of $g$ we also have:
  \begin{equation*}
    \quad \vdash_{PA} \mathcal{D}(\overline{m}, \overline{g})\quad \quad   \text{(\ref{proposition_godel_sentence}.4, \ref{lemma_diagonal_relation_sign}.1, modus ponens)} \tag{\ref{proposition_godel_sentence}.8}
  \end{equation*}
  We now proceed to a proof by cases.
  \begin{description}
  \item[Case $\vdash_{PA} \mathcal{G} \Rightarrow (\forall x_2) (\lnot B(x_2, \overline{g}))$]
    \begin{equation*}
      (\forall z) \{\mathcal{D}'(\overline{m}, z) \Leftrightarrow (\forall x_2) [\lnot B'(x_2, z)]\} \quad \quad \quad \quad \quad \quad \quad \quad \quad  \tag{1} \text{hypothesis }
    \end{equation*}
    \begin{multline*}
      \quad \vdash_{PA} \{\mathcal{D}'(\overline{m}, \overline{g}) \Leftrightarrow (\forall x_2) [\lnot B'(x_2, \overline{g})]\} \tag{2}\ \text{(1), $A4$ / universal instantiation}
    \end{multline*}
    \begin{multline*}
      \quad \vdash_{PA} (\forall x_2) [\lnot B'(x_2, \overline{g})] \tag{3}\ \\
      \text{(2), \ref{proposition_godel_sentence}.8, biconditional elimination}
    \end{multline*}
    \begin{multline*}
      \quad (\forall z) \{\mathcal{D}'(\overline{m}, z) \Leftrightarrow (\forall x_2) [\lnot B'(x_2, z)]\} \\
      \vdash_{PA} (\forall x_2) [\lnot B'(x_2, \overline{g})] \tag{4}\ \text{(1)-(3), transitivity of conditional proof}
    \end{multline*}
    \begin{multline*}
      \quad \vdash_{PA} (\forall z) \{\mathcal{D}'(\overline{m}, z) \Leftrightarrow (\forall x_2) [\lnot B'(x_2, z)]\} \\
      \Rightarrow (\forall x_2) [\lnot B'(x_2, \overline{g})] \tag{5}\ \text{(4), deduction theorem}
    \end{multline*}
  \item[Case $\vdash_{PA} (\forall x_2)(\lnot B(x_2, \overline{g})) \Rightarrow \mathcal{G} $]% (\forall x_2) $]
    \begin{equation*}
      (\forall x_2)(\lnot B'(x_2, \overline{g})) \quad \quad \quad \quad \quad \quad \quad \quad \quad \quad \quad \quad \quad \tag{7} \text{hypothesis }
    \end{equation*}
    \begin{equation*}
      \mathcal{D}'(\overline{m}, z) \quad \quad \quad \quad  \quad \quad \quad \quad  \quad \quad \quad \quad \quad \quad \quad \quad \tag{8} \text{hypothesis }
    \end{equation*}
    \begin{equation*}
      \vdash_{PA} (\exists_1 z) \mathcal{D}'(\overline{m}, z) \quad \quad \quad \quad  \quad \quad \quad \quad  \quad \quad \quad \quad \quad \quad \tag{9} \text{\ref{lemma_diagonal_relation_sign}.3}
    \end{equation*}
    \begin{flalign*}
      \vdash_{PA} (\exists z) \mathcal{D}'(\overline{m}, z)\ \land (\forall x)(\forall y) [\mathcal{D}'(\overline{m}, x) \land  
        \mathcal{D}'(\overline{m}, y) \Rightarrow x = y]\  \tag{10} \text{(9), (D5), }\\ %\quad \quad \quad \quad
      \text{equivalence by definition, conjunction elimination, modus ponens} 
    \end{flalign*}
    \begin{flalign*}
      \vdash_{PA} (\forall x)(\forall y) [\mathcal{D}'(\overline{m}, x) \land  
        \mathcal{D}'(\overline{m}, y) \Rightarrow x = y]  \tag{11} \\ %\quad \quad \quad \quad
      \text{(10), conjunction elimination}
    \end{flalign*}
    \begin{flalign*}
      \vdash_{PA} (\forall y) [\mathcal{D}'(\overline{m}, \overline{g}) \land  
        \mathcal{D}'(\overline{m}, y) \Rightarrow \overline{g} = y]  \tag{12} \\ %\quad \quad \quad \quad
      \text{(11), $A4$ / universal instantiation}
    \end{flalign*}
    \begin{flalign*}
      \vdash_{PA} [\mathcal{D}'(\overline{m}, \overline{g}) \land  
        \mathcal{D}'(\overline{m}, z) \Rightarrow \overline{g} = z]  \tag{13} \\ %\quad \quad \quad \quad
      \text{(12), $A4$ / universal instantiation}
    \end{flalign*}
    \begin{flalign*}\label{equation_godel_sentence_onlyif_conjunction_introduction}
      \vdash_{PA} \mathcal{D}'(\overline{m}, \overline{g}) \land  
      \mathcal{D}'(\overline{m}, z)\  \tag{14} 
      \text{(\ref{proposition_godel_sentence}.8, (8), conjunction introduction)}
    \end{flalign*}
    \begin{equation*}
      \vdash_{PA} (\overline{g} = z) \quad \quad  \quad \quad \quad \quad \quad \quad \tag{15} \text{(13), (14), modus ponens}
    \end{equation*}
    \begin{equation*}
      \vdash_{PA} (\forall x_2)[\lnot B'(x_2, z)] \quad \quad \quad \quad \quad \quad \tag{16} \text{(7), (15), substitutivity of equality}
    \end{equation*}
    \begin{multline*}
      \quad (\forall x_2)(\lnot B(x_2, \overline{g})), \mathcal{D}'(\overline{m}, z) \vdash_{PA} (\forall x_2)[\lnot B'(x_2, z)] \tag{17}\ \\
      \text{(7)-(16), transitivity of conditional proof} 
    \end{multline*}
    \begin{multline*}
      \quad (\forall x_2)(\lnot B(x_2, \overline{g})) \vdash_{PA} \mathcal{D}'(\overline{m}, z) \Rightarrow (\forall x_2)[\lnot B'(x_2, z)] \tag{18}\ \\
      \text{(17), deduction theorem} 
    \end{multline*}
    \begin{multline*}
      \quad (\forall x_2)(\lnot B(x_2, \overline{g})) \vdash_{PA} (\forall z) \{\mathcal{D}'(\overline{m}, z) \Rightarrow (\forall x_2)[\lnot B'(x_2, z)]\} \tag{19}\ \\
      \text{(18), generalization} 
    \end{multline*}
    \begin{multline*}
      \quad \vdash_{PA} (\forall x_2)(\lnot B(x_2, \overline{g})) \Rightarrow (\forall z) \{\mathcal{D}'(\overline{m}, z) \Rightarrow (\forall x_2)[\lnot B(x_2, z)]\} \tag{20}\ \\
      \text{(19), deduction theorem} 
    \end{multline*}
    At (\ref{equation_godel_sentence_onlyif_conjunction_introduction}) and elsewhere, use is made of the fact that, for example, 
    $\mathrm{Subst}\ \left[ \mathcal{D}'\ \genfrac{}{}{0pt}{0}{x_1}{\overline{m}} \genfrac{}{}{0pt}{0}{z}{\overline{g}} \right]$
    and 
    $\mathrm{Subst}\ \left[ \mathcal{D}\ \genfrac{}{}{0pt}{0}{x_1}{\overline{m}} \genfrac{}{}{0pt}{0}{x_2}{\overline{g}} \right]$
    are the same formula.
  \end{description}

\end{proof}
Whilst Equation \ref{proposition_godel_sentence}.1 is generally used to refer to the G\"odel sentence for $PA$,  $\mathcal{G}$, for our purposes
it will be useful to note the following $PA$ theorem concerning this sentence:
\begin{Corollary} \label{corollary_godel_sentence}
  For the $PA$ G\"odel-sentence $\mathcal{G}$, with G\"odel number $g$, defined at Proposition \ref{proposition_godel_sentence},
  we also have:
  \begin{equation*}
    \vdash_{PA} (\lnot \mathcal{G}) \Leftrightarrow (\exists x_2) [B(x_2, \overline{g})] \tag{\ref{corollary_godel_sentence}.1}
  \end{equation*}
\end{Corollary}
\begin{proof}
  \begin{quote}
    \begin{multline*}
      \quad (1)\ \vdash_{PA} (\lnot \mathcal{G}) \Leftrightarrow (\lnot \{(\forall x_2) [\lnot B(x_2, \overline{g})]\}) \\
      (\ref{proposition_godel_sentence}.1),\ \text{biconditional negation}
    \end{multline*}
    \begin{multline*}
      \quad (2)\ \vdash_{PA} ( (\exists x_2) [B(x_2, \overline{g})] ) \Leftrightarrow (\lnot \{(\forall x_2) [\lnot B(x_2, \overline{g})]\})\ \\
      \text{definition of } (\exists x) \mathcal{B}(x)\text{,} \text{(\cite{mendelson2015}: 48), equivalence by definition}
    \end{multline*}
    \begin{multline*}
      \quad (3)\ \vdash_{PA} (\lnot \mathcal{G}) \Leftrightarrow (\exists x_2) [B(x_2, \overline{g})]  \\
      \text{ (1), (2),} \text{Equivalence theorem (\cite{mendelson2015}: Proposition 2.9(a))}
    \end{multline*}
  \end{quote}
\end{proof}
With the above material at hand some new results about the arithmetization of the syntax of first-order number theory may be established.
\subsection{A novel demonstration}
For the last proposition required for the demonstration I have not as yet
identified a purely syntactic demonstration in the existing literature. I will
firstly consider such a proof before considering whether an existing
demonstration from the literature that uses semantic methods supports this proposition.
\begin{Proposition} \label{proposition_proof_inference}
  If $\mathcal{F}_y$ is any $PA$ formula with G\"odel number $y$, then
  Propositions \ref{proposition_arithmetized_proof} and \ref{proposition_proof_relation_sign} imply,
  on the assumption that $PA$
  is $\omega$ consistent (or consistent, respectively), that the following propositions both hold:
  \begin{equation*}
    \{\vdash_{PA}  (\exists x_2) [B(x_2, \overline{y})]\} \rightarrow\ \vdash_{PA} \mathcal{F}_y \tag{\ref{proposition_proof_inference}.1}
  \end{equation*}
  \begin{equation*}
    \{\vdash_{PA}  (\forall x_2) [\lnot B(x_2, \overline{y})]\} \rightarrow\ \nvdash_{PA} \mathcal{F}_y \tag{\ref{proposition_proof_inference}.2}
  \end{equation*}
\end{Proposition}
\begin{proof} Since Proposition \ref{proposition_proof_inference} is metatheoretical in nature the proof is naturally metatheoretical.
  The proof by contradiction in the first case (\ref{proposition_proof_inference}.1) assumes that that $PA$ is $\omega$ consistent and in the
  second case (\ref{proposition_proof_inference}.2) that $PA$ is simply consistent. (The proof in the first case parallels G\"odel's argument
  that $\mathrm{Neg} [17\ \mathrm{Gen}\ r]$ is not $\kappa$-\textsc{provable}, and in the second case parallels G\"odel's argument
  that $17\ \mathrm{Gen}\ r$ is not $\kappa$-\textsc{provable}.)
  \begin{description}
  \item[Proof of (\ref{proposition_proof_inference}.1)] Suppose that Equation \ref{proposition_proof_inference}.1 fails.
    By hypothesis, \\
    $\vdash_{PA} (\exists x_2) [B(x_2, \overline{y})]$ is true and
    yet 
    $\vdash_{PA} \mathcal{F}_y$ is false. But if $\nvdash_{PA} \mathcal{F}_y$ thus holds, then there is no $PA$ proof of $\mathcal{F}_y$
    and thus, by the hypothesis that the relation $x\ \mathrm{B}\ y$
    corresponds to the $PA$ proof relation, $\overline{n\ \mathrm{B}\ y}$ holds for every natural number $n$. Hence furthermore,
    by Equation \ref{proposition_proof_relation_sign}.2, $(n)\vdash_{PA} \lnot [B(\overline{n}, \overline{y})]$ holds.
    The hypothesis that $\vdash_{PA} (\exists x_2) [B(x_2, \overline{y})]$
    holds thus contradicts the assumption that $PA$ is $\omega$ consistent.

  \item[Proof of (\ref{proposition_proof_inference}.2)] Suppose that Equation \ref{proposition_proof_inference}.2 fails.
    By hypothesis, \\
    $\vdash_{PA} (\forall x_2) [\lnot B(x_2, \overline{y})]$ is true and
    yet $\nvdash_{PA} \mathcal{F}_y$ is false. But if $\vdash_{PA} \mathcal{F}_y$ holds, then there is a $PA$ proof of $\mathcal{F}_y$.
    Thus, by the hypothesis that the relation $x\ \mathrm{B}\ y$
    corresponds to the $PA$ proof relation, there exists a natural number $n$ such that $n\ \mathrm{B}\ y$ holds. Hence,
    by Equation \ref{proposition_proof_relation_sign}.1, $\vdash_{PA} B(\overline{n}, \overline{y})$ also holds.
    Since $\vdash_{PA} (\forall x_2) [\lnot B(x_2, \overline{y})]$ implies by universal instantiation (\cite{mendelson2015} A4) and modus ponens
    that $\vdash_{PA} [\lnot B(\overline{n}, \overline{y})]$ also holds,
    the hypothesis for this case is incompatible with the assumption that $PA$ is consistent.
  \end{description}
  
  To formalise the demonstration let's consider the arithmetic image of (\ref{proposition_proof_inference}.1-2),
  using some further symbolism and material from \cite{godel1931}.
  \begin{enumerate}
  \item Let $v_1$ and $v_2$ be the G\"odel numbers of the first and second variables of $PA$ respectively;
  \item Let $b$  the G\"odel number of the $PA$ relation sign $B(x, y)$ mentioned at Proposition \ref{proposition_proof_relation_sign};
  \item Where $f$ is the G\"odel number of a $PA$ formula $\mathcal{F}$:
    \begin{enumerate}
    \item Let $ \mathrm{Neg} (f)$ be the G\"odel number of $(\lnot \mathcal{F})$;
    \item Let $ v_i\ \mathrm{Gen}\ f$ the G\"odel number of $(\forall x_i) \mathcal{F}$;
    \item Let $ v_i\ \mathrm{Ex}\ f$ the G\"odel number of $(\exists x_i) \mathcal{F}$.
    \end{enumerate}
  \end{enumerate}
  If we then consider the arithmetic image of $PA$, Equations \ref{proposition_proof_inference}.1-2 correspond to the following arithmetic propositions:
  \begin{equation*}
    (y)(\{ (Ez) [z\ B\ (v_2\ \mathrm{Ex}\ \{ \mathrm{Sb}\ \left[ b\ \genfrac{}{}{0pt}{0}{v_1}{Z(y)} \right] \} )]\} \rightarrow\ \{ (Ez) [z\ B\ y]\}) \tag{\ref{proposition_proof_inference}.1.1}
  \end{equation*}
  \begin{equation*}
    (y)(\{ (Ez) [z\ B\ (v_2\ \mathrm{Gen}\ \{ \mathrm{Neg}\ [ \mathrm{Sb}\ \left[ b\ \genfrac{}{}{0pt}{0}{v_1}{Z(y)} \right] ] \} )]\} \rightarrow\ \{ \overline{(Ez)} [z\ B\ y]\}) \tag{\ref{proposition_proof_inference}.2.1}
  \end{equation*}

  \paragraph{Proof of (\ref{proposition_proof_inference}.1.1)}\
  The falsity of Equation \ref{proposition_proof_inference}.1.1 implies:
  \begin{flalign*}
    \quad (1)\ (Ey)(\{ (Ez) [z\ B\ (v_2\ \mathrm{Ex}\ \{ \mathrm{Sb}\ \left[ b\ \genfrac{}{}{0pt}{0}{v_1}{Z(y)} \right] \} )]\}\ \&\ \{ \overline{(Ez)} [z\ B\ y]\})\ && \\
    \text{Since, }  \{\overline{(x)}[F(x) \rightarrow G(x)] \equiv (Ex) [F(x)\ \&\ \overline{G}(x)]\}
  \end{flalign*}
  \begin{flalign*}
    \quad (2)\ (\{ (Ez) [z\ B\ (v_2\ \mathrm{Ex}\ \{ \mathrm{Sb}\ \left[ b\ \genfrac{}{}{0pt}{0}{v_1}{Z(c)} \right] \} )]\}\ \&\ \{ \overline{(Ez)} [z\ B\ c]\})\ && \\
    \text{(1), Existential instantiation, c is new}
  \end{flalign*}
  \begin{flalign*}
    \quad (3)\ (Ez) [z\ B\ (v_2\ \mathrm{Ex}\ \{ \mathrm{Sb}\ \left[ b\ \genfrac{}{}{0pt}{0}{v_1}{Z(c)} \right] \} )]\ &&  \text{(2), conjunction elimination}
  \end{flalign*}
  \begin{flalign*}
    \quad (4)\ [p\ B\ (v_2\ \mathrm{Ex}\ \{ \mathrm{Sb}\ \left[ b\ \genfrac{}{}{0pt}{0}{v_1}{Z(c)} \right] \} )]\  && \text{(3), existential instantiation, p is new}
  \end{flalign*}
  \begin{flalign*}
    \quad (5)\ \overline{(Ez)} [z\ B\ c]\ && \text{(2), conjunction elimination}
  \end{flalign*}
  \begin{flalign*}
    \quad (6)\ (z) \overline{[z\ B\ c]}\ && \text{(5), definition of}\ (Ez)
  \end{flalign*}
  \begin{flalign*}
    \quad (7)\ (z)\{ (Ew) [w\ B\ ( \mathrm{Neg}\ \{\mathrm{Sb}\ \left[ b\ \genfrac{}{}{0pt}{0}{v_1}{Z(c)} \genfrac{}{}{0pt}{0}{v_2}{Z(z)} \right] \} )] \} && \text{(6), (\ref{proposition_proof_relation_sign}.2), modus ponens}  
  \end{flalign*}
  (4) and (7) contradict the assumption that $PA$ is $\omega$ consistent. Since the falsity of Equation \ref{proposition_proof_inference}.1.1 implies a contradiction,
  the truth of Equation \ref{proposition_proof_inference}.1.1 may be inferred. \\
  
  \paragraph{Proof of (\ref{proposition_proof_inference}.2.1)}\
  The falsity of Equation \ref{proposition_proof_inference}.2.1 implies:
  \begin{flalign*}
    \quad (1)\ (Ey)(\{ (Ez) [z\ B\ (v_2\ \mathrm{Gen}\ \{ \mathrm{Neg}\ [\mathrm{Sb}\ \left[ b\ \genfrac{}{}{0pt}{0}{v_1}{Z(y)} \right]] \} )] \}\ \&\ \{ (Ez) [z\ B\ y]\})\ && \\
    \text{Since, } \{\overline{(x)}[F(x) \rightarrow G(x)] \equiv (Ex) [F(x)\ \&\ \overline{G}(x)]\}
  \end{flalign*}
  \begin{flalign*}
    \quad (2)\ \{ (Ez) [z\ B\ (v_2\ \mathrm{Gen}\ \{ \mathrm{Neg}\ [\mathrm{Sb}\ \left[ b\ \genfrac{}{}{0pt}{0}{v_1}{Z(c)} \right]] \} )] \}\ \&\ \{ (Ez) [z\ B\ c]\}\ && \\
    \text{(1), Existential instantiation, c is new}
  \end{flalign*}
  \begin{flalign*}
    \quad (3)\ (Ez) [z\ B\ (v_2\ \mathrm{Gen}\ \{ \mathrm{Neg}\ [\mathrm{Sb}\ \left[ b\ \genfrac{}{}{0pt}{0}{v_1}{Z(c)} \right]]\} ) ]\ &&  \text{(2), conjunction elimination}
  \end{flalign*}
  \begin{flalign*}
    \quad (4)\ (Ez) [z\ B\ c]\ && \text{(2), conjunction elimination}
  \end{flalign*}
  \begin{flalign*}
    \quad (5)\ [p\ B\ c]\  && \text{(4), existential instantiation, p is new}
  \end{flalign*}
  \begin{flalign*}
    \quad (6)\ (Ew) [w\ B\ (\mathrm{Sb}\ \left[ b\ \genfrac{}{}{0pt}{0}{v_1}{Z(c)} \genfrac{}{}{0pt}{0}{v_2}{Z(p)} \right])] \} && \text{(5), (\ref{proposition_proof_relation_sign}.2), modus ponens}  
  \end{flalign*}
  \begin{flalign*}
    \quad (7)\ (Eu) [u\ B\ (\mathrm{Neg}\ [\mathrm{Sb}\ \left[ b\ \genfrac{}{}{0pt}{0}{v_1}{Z(c)} \genfrac{}{}{0pt}{0}{v_2}{Z(p)} \right]]\ &&  \text{(3), Particularisation rule $A4$ (\cite{mendelson2015}: 2.5.1) }
  \end{flalign*}
  (6) and (7) contradict the assumption that $PA$ is consistent.  Since the falsity of Equation \ref{proposition_proof_inference}.2.1 implies a contradiction,
  the truth of Equation \ref{proposition_proof_inference}.2.1 may be inferred. \\
\end{proof}

I will now look at how
Proposition \ref{proposition_proof_inference}
may be used to establish
the main result before considering (in \S \ref{subsection_derive_established_results}) whether
the former proposition can be established
using existing results from the literature.\footnote{The digression
into existing related results may not be of interest to a reader who accepts the
above proof of Proposition \ref{proposition_proof_inference}.}%\newpage
\begin{Proposition}\label{propostion_theorems_not_defined}
  The class of $PA$ theorems is not well defined.
\end{Proposition}
\begin{proof}
  \begin{quote}
    \begin{flalign*}
      \quad (1)\ (\lnot \mathcal{G})  && \text{Hypothesis}
    \end{flalign*}
    \begin{flalign*}
      \quad (2)\  (\lnot \mathcal{G}) \vdash_S (\lnot \mathcal{G}) && \text{(1),}\ \text{Conditional identity}
    \end{flalign*}
    \begin{flalign*}
      \quad (3)\ \vdash_{PA} (\exists x_2) B(x_2, \overline{g}) && \text{(2),}\ (\ref{corollary_godel_sentence}.1),\ \text{biconditional elimination}
    \end{flalign*}
    \begin{flalign*}
      \quad (4)\ \vdash_{PA} \mathcal{G} && \text{(3),}\ (\ref{proposition_proof_inference}.1), \text{modus ponens}
    \end{flalign*}
    \begin{flalign*}
      \quad (5)\ (\lnot \mathcal{G})  \vdash_{PA} \mathcal{G} && \text{transitivity of conditional proof}\ (1)-(4)  
    \end{flalign*}
    \begin{flalign*}
      \quad (6)\ \vdash_{PA} (\lnot \mathcal{G}) \Rightarrow \mathcal{G} && (5), \text{deduction theorem}  
    \end{flalign*}
    \begin{flalign*}
      \quad (7)\ \vdash_{PA} [(\lnot \mathcal{G}) \Rightarrow \mathcal{G}] \Rightarrow\ \mathcal{G} && \text{tautology theorem}  
    \end{flalign*}
    \begin{flalign*}
      \quad (8)\ \vdash_{PA} \mathcal{G} && (6), (7), \text{modus ponens}  
    \end{flalign*}
    \begin{flalign*}
      \quad (9)\ \vdash_{PA} (\forall x_2) [\lnot B(x_2, \overline{g})] && (8), (\ref{proposition_godel_sentence}.1), \text{biconditional elimination}  
    \end{flalign*}
    \begin{flalign*}
      \quad (10)\ \nvdash_{PA} \mathcal{G} && (9), (\ref{proposition_proof_inference}.2), \text{modus ponens}  
    \end{flalign*}
  \end{quote}
  Since the theoremhood
  of $\mathcal{G}$ is both proved at (8) and refuted at (10)
  the class of $PA$ theorems is clearly not well defined.
\end{proof}
The reader who is not convinced of the
above proof of Proposition \ref{proposition_proof_inference}
may wish to consider whether this Proposition
can be otherwise derived from results already
established in the literature. The reader who accepts the
above proof may skip \S \ref{subsection_derive_established_results} without loss.

\subsection{A derivation of Proposition \ref{proposition_proof_inference} from established results}\label{subsection_derive_established_results}
To obtain a proof of Proposition \ref{proposition_proof_inference} from the existing literature,
we can use \'Swierczkowski's "Proof formalization condition",
a metatheorem he establishes for hereditarily finite set theory ($HF$):
\begin{quote}
  \textsc{Proposition 4.4} (Proof formalization condition). For every formula $\varphi$,
  \begin{flalign*}
    \quad \quad  \quad \quad \quad \vdash \varphi\ \text{iff} \vdash \mathrm{Pf}(\ulcorner \varphi \urcorner) \text{. (\cite{swierczkowski2003}: 18)} &&
  \end{flalign*}
\end{quote}
As my aim is only to indicate an alternative method of establishing Proposition \ref{proposition_proof_inference},
I simply sketch an alternative demonstration that uses only already accepted results.

Focusing firstly on $HF$ itself we can note the following.
\begin{Proposition} \label{proposition_proof_formalization_implies_proof_inference_hf}
  The Proof formalization condition (Proposition 4.4) implies,
  if $HF$ is assumed consistent, that a metatheorem
  corresponding to Proposition \ref{proposition_proof_inference} applies to $HF$.
\end{Proposition}
\begin{proof}
  If we decompose the Proof formalization biconditional into its two conjuncts we obtain:
  \begin{flalign*}
    \quad \quad \quad \quad \quad  \{\vdash \mathrm{Pf}(\ulcorner \varphi \urcorner)\} \rightarrow \{\vdash \varphi\ \} && \tag{4.4.1}
  \end{flalign*}
  \begin{flalign*}
    \quad \quad \quad \quad \quad \{\vdash \varphi\ \} \rightarrow \{\vdash \mathrm{Pf}(\ulcorner \varphi \urcorner)\} && \tag{4.4.2}
  \end{flalign*}
  Equation 4.4.1 corresponds to Equation \ref{proposition_proof_inference}.1.
  The Equation corresponding to Equation \ref{proposition_proof_inference}.2 is this:
  \begin{flalign*}
    \quad \quad \quad \quad \quad  \{\vdash \lnot \mathrm{Pf}(\ulcorner \varphi \urcorner)\} \rightarrow \{\nvdash \varphi\ \} && \tag{4.4.3}
  \end{flalign*}
  It may be seen that (4.4.2) implies
  that (4.4.3) holds, given the assumption that $HF$ is consistent,
  by the following variation of the demonstration given above.
  Suppose to derive a contradiction that (4.4.3) fails.
  By hypothesis, 
  $\{\vdash \lnot \mathrm{Pf}(\ulcorner \varphi \urcorner)\}$ is true and
  yet $\{\nvdash \varphi\ \}$ is false. But if  $\{\vdash \varphi\ \}$ holds,
  then by (4.4.2) we have $\{\vdash \mathrm{Pf}(\ulcorner \varphi \urcorner)\}$,
  which, in conjunction with the first hypothesis,
  contradicts the assumption that $HF$ is consistent.
  Hence, Equation 4.4.3 holds.
  Thus, given the Proof formalization condition and
  the assumption that $HF$ is consistent,
  a metatheorem  corresponding to Proposition \ref{proposition_proof_inference} applies to $HF$
  also.
\end{proof}
In a previous context the concept of an
the concept of an an extension by definitions was introduced
(Definition \ref{definition_extension_by_definitions}).
\'Swierczkowski's Theorem 10.1 establishes that "each of the theories HF and
PA is an extension by definitions of the other" (\cite{swierczkowski2003}:
36). Theorem 10.1 essentially
defines a formal theory $\mathcal{T}$ that is a common extension
by definitions of both $PA$ and $HF$. The details are
omitted here since all that is required is the existence of Theorem 10.1
as an accepted result supporting the cited claim.
From Theorem 10.1 we may infer the following:
\begin{Corollary} \label{corollary_proof_formalization_implies_proof_inference}
  The Proof formalization condition (Proposition 4.4),
  in conjunction with Theorem 10.1 (\cite{swierczkowski2003}: 37)
  and the assumption that $HF$ consistent, 
  implies that Equations \ref{proposition_proof_inference}.1-2 hold
  for an arbitrary $PA$ formula  $\mathcal{F}_y$ with G\"odel number $y$.
\end{Corollary}
\begin{proof}
  Since $\mathcal{T}$ is an extension
  by definitions of both $PA$ and $HF$,
  we have, by Definition \ref{definition_extension_by_definitions}:
  \begin{Lemma} \label{lemma_hf_pa_common_extension}
    \begin{enumerate}
    \item \label{item_hf_pa_common_extension_hf} For each formula $\mathcal{A}$ of $\mathcal{T}$,
      there exists a corresponding
      formula  $\mathcal{A}^{*}$ of $HF$ such that:
      \begin{equation*}
        \vdash_{HF} \mathcal{A}^{*} \leftrightarrow\ \vdash_{\mathcal{T}} \mathcal{A}
      \end{equation*}
    \item \label{item_hf_pa_common_extension_pa} For each formula $\mathcal{A}$ of $\mathcal{T}$,
      there exists a corresponding
      formula  $\mathcal{A}^{*}$ of $PA$ such that:
      \begin{equation*}
        \vdash_{PA} \mathcal{A}^{*} \leftrightarrow\ \vdash_{\mathcal{T}} \mathcal{A}
      \end{equation*}
  \end{enumerate}
\end{Lemma}
From Item \ref{item_hf_pa_common_extension_hf} it follows immediately
that the  Proof formalization condition (Proposition 4.4),
with appropriate changes, applies to $\mathcal{T}$;
this, in conjunction with From Item \ref{item_hf_pa_common_extension_pa}
implies that the  Proof formalization condition (Proposition 4.4),
with appropriate changes, applies to $PA$.
\end{proof}
At this point a minor qualification should be noted.
\'Swierczkowski's proof of the Proof formalization condition
uses a semantic claim, namely that for every $HF$ sentence $\phi$ of a certain form
truth in a standard model $\mathfrak{S}$ of the theory implies provability in $HF$:
$ \{\mathfrak{S}\vDash_{HF}  \phi \} \rightarrow\ \vdash_{HF} \phi$.

Thus if \'Swierczkowski's proof of the Proof formalization condition is relied
upon, rather than the above syntactic proof, a purely syntactic proof is not available.

Now a second minor qualification should be noted. Proposition \ref{proposition_proof_inference}
bundles together two derived inference rules (\ref{proposition_proof_inference}.1
and \ref{proposition_proof_inference}.2). The first of these follows immediately
from \'Swierczkowski's Proof formalization condition (Proposition 4.4) - given
an auxiliary, uncontroversial proposition concerning the relationship
between $PA$ and $HF$. To derive the second
Proposition \ref{proposition_proof_inference} inference
rule (Equation \ref{proposition_proof_inference}.2) from the Proof formalization condition,
the very brief deduction,
presented in the proof of Proposition \ref{proposition_proof_formalization_implies_proof_inference_hf}, 
is used; this deduction has not as yet been independently confirmed.

Suppose, for the sake of argument,
that the reader does not accept the very brief deduction
presented in the proof of Proposition \ref{proposition_proof_formalization_implies_proof_inference_hf}.
In the following subsection I show that
\'Swierczkowski's Proof formalization condition (Proposition 4.4)
nevertheless still yields Proposition \ref{propostion_theorems_not_defined}.
\subsubsection{The derivation with (\ref{proposition_proof_inference}.2) removed}
For emphasis, the aim of the following discussion is simply to clarify
the impact on the above proof of the removal of Equation \ref{proposition_proof_inference}.2.
There is no suggestion that there is any problem with either this Equation or
the two proofs of it presented above; neither proof however has been independently confirmed.
If the reader is satisfied with either demonstration then they may skip the
remainder of this subsection without loss.

From what
is stated above, with everything that depends upon Equation \ref{proposition_proof_inference}.2
removed, we may still infer the following:
\begin{Corollary} \label{corollary_proof_formalization_implies_proof_inference_1}
  The Proof formalization condition (Proposition 4.4),
  in conjunction with Theorem 10.1 (\cite{swierczkowski2003}: 37)
  implies that Equation \ref{proposition_proof_inference}.1 holds
  for an arbitrary $PA$ formula  $\mathcal{F}_y$ with G\"odel number $y$.
\end{Corollary}
Proceeding then on the basis of the more restricted
Corollary \ref{corollary_proof_formalization_implies_proof_inference_1}, we have: 
\begin{Corollary}\label{corollary_theorems_not_defined}
  Proposition \ref{propostion_theorems_not_defined} 
  can be obtained assuming only results established by accepted proofs.
\end{Corollary}
\begin{proof}
  Steps (1)-(8) as set out in the above proof of Proposition \ref{propostion_theorems_not_defined}
  are the same. From this point however we may adapt the standard
  proof of incompleteness to our purpose (\cite{mendelson2015}: Proposition 3.37,
  close paraphrase follows). From $\vdash_{PA} \mathcal{G}$ (8), by Rule C in the metatheory,
  let $r$ be the G\"odel number of this proof; that is, from $(Ez) [z\ \mathrm{B}\ g]$
  we chose $r$ so that $r\ \mathrm{B}\ g$ holds. Hence $\vdash_{PA} B(\overline{r}, \overline{g})$
  follows by Equation \ref{proposition_proof_relation_sign}.1.
  The alternative deduction then continues as follows.
  \begin{quote}
    \begin{flalign*}
      \quad (9')\ \vdash_{PA} B(\overline{r}, \overline{g})  && \text{(8), Rule C in metatheory, (\ref{proposition_proof_relation_sign}.1})
    \end{flalign*}
    \begin{flalign*}
      \quad (10')\ \vdash_{PA} (\forall x_2)[\lnot B(x_2, \overline{g})] && \text{(8),}\ (\ref{proposition_godel_sentence}.1),\ \text{biconditional elimination}
    \end{flalign*}
    \begin{flalign*}
      \quad (11')\ \vdash_{PA} [\lnot B(\overline{r}, \overline{g})] && \text{(12), Particularisation rule $A4$ (\cite{mendelson2015}: 2.5.1)}
    \end{flalign*}
  \end{quote}
\end{proof}
(9') and (11') imply, strictly speaking, that $S$ is formally inconsistent.
If we adopt a classical approach to logic, and assume that the informal
arithmetic of natural numbers is consistent, then we may take this result as
an alternative confirmation of
a failure of the metamathematical definition of $PA$'s proof relation
and hence an alternative confirmation of Proposition \ref{propostion_theorems_not_defined}.

\subsection{Discussion}\label{subsection_discussion}

The breakdown in the metamathematical definition of $PA$'s proof relation points to a failure in the
so-called simple theory of types relied upon for the prevention of paradox. Simple type theory, as opposed to the ramified type theory of the first edition
of \emph{Principia Mathematica} \cite{pm1910v1}, distinguishes semantic paradoxes, involving such notions such as truth, from paradoxes which
involve only notions involved in ordinary mathematics such as set theory (\cite{mendelson2015}: Introduction, cf. \cite{ramsey1925}: 20). The above paradox falls into the latter camp,
since the metamathematical notions involved, of proof etc., do not involve semantic notions.

Metamathematical presentations of type theory tend to obscure an important distinction between \emph{Principia}'s type theory, on
the one hand, and the metamathematical approach, exemplified for example in G\"odel's system $P$, on the other.

On the metamathematical approach, types are assigned to uninterpreted
symbols of the formal system alone. To use an example from Mendelson's first-order number theory $S$ (\cite{mendelson2015}: Chapter 3):
Under the standard interpretation $\mathfrak{M}$ the sentence ``$\overline{1} = \overline{1}$'' is associated with the (complex) of the number one (whatever that may be,
in the domain of interpretation) standing in the relation of equality to itself. The type of ``being an individual'' for $S$, on the metamathematical approach, is a property of the
relevant symbols (symbol-sequences) themselves - the symbol for zero, as well as symbols for individual variables, and, indirectly in a sense, the
symbol sequences for other numerals and terms. Under a non-standard interpretation of $S$, the objects assigned to these symbols,
contained in the domain of interpretation, need not bear any relation to the objects that are members of the domain of the standard
interpretation; these objects might be classes of numbers, or sequences of numbers, or uninterpreted symbols, or just about anything.

On \emph{Principia}'s approach, types are assigned to the
\emph{interpreted} symbols of both the formal system being defined and the informal language used to thus define \emph{Principia}'s formal system.
In \emph{Principia}'s type theory, when members of the type of
individuals are spoken of, in a given context, it is the symbols for individuals in association with the objects assigned to these under
\emph{Principia}'s fixed interpretation that are mentioned. The individuals themselves, which are mentioned by using these symbols, are essentially, anything
that is neither a proposition nor a propositional function (\cite{boyce2025}).

This distinction is crucially important, since G\"odel's definition of the syntax of his system $P$, and a fortiori his arithmetization of this syntax,
requires the metamathematical approach to type theory, whilst \emph{Principia}'s
type theory blocks this arithmetization. The arithmetization of syntax requires, for example, that there exists a relation of ``immediate consequence'' that holds
between all formulae $a_1$, $a_2$,  $a_3$ of the system such that $a_3$ may be inferred from $a_1$ / $a_2$ according to the rules of inference of the system:
\begin{quote} The class of \emph{provable formulas} is defined to be the
  smallest class of formulas that contains the axioms and is closed
  under the relation "immediate consequence". (\cite{godel1931}: 155-7)
\end{quote}
G\"odel's seemingly innocuous assumption that a well-defined class of provable formulae, and associated notion of formal proof, is determined
in the above-quoted manner is directly refuted by the above result.

By contrast, \emph{Principia}'s axioms,
propositions and propositional functions, as indicated in the phrasing of it's Primitive propositions (rules of inference / formation rules) are \emph{interpreted} and \emph{not}
all of the same type (\cite{boyce2024}). Consider for example \emph{Principia}'s statement of modus ponens for elementary propositions:
\begin{quote}
  $\pmast 1\pmcdot1$ Anything implied by a true elementary proposition is true. \pmpp. (\cite{pm1910v1}: 98)
\end{quote}
Whilst the principle applies to all elementary propositions, it cannot meaningfully be applied to propositions of a higher type, for which purpose
a different Primitive Proposition is required. Whitehead and Russell's formulation of the principle looks quite strange to contemporary
eyes accustomed to the metamathematical approach. Many subtleties are involved which fall outside the focus of discussion here, including the fact
that: since the (interpreted) natural language in which Whitehead and Russell define and reason about their formal system is also typed,
it is in scope for the formal system defined (taking into account the types of the propositions / propositional functions involved).

The important point for present purposes is that by \emph{Principia}'s theory of logical types,
\emph{Principia}'s formulae and theorems cannot all be members of a single class. The formulae and theorems are not all of the same type,
and, according to \emph{Principia}'s theory of classes,
members of a class must all be of the same type. On this view, G\"odel's assumption that one can define a single relation of
"immediate consequence" - which can then be arithmetized - is false. \emph{Principia}'s arithmetic thus avoids the above paradox.
It is thus somewhat ironic that this feature of \emph{Principia},
in virtue of which the above paradox is avoided, is also a feature which has attracted rather harsh criticism from advocates of metamathematics (\cite{godel1944}: 120).

As \emph{Principia}'s type theory prevents the indirect self-reference that underpins the above paradox we may say that by \emph{Principia}'s lights the prohibited propositions /
propositional functions are meaningless. \emph{Principia}'s type theory thus provides a uniform treatment of both L\"ob's theorem and Curry's paradox (\cite{lob1955}). The formulae
in scope for L\"ob's theorem, like the G\"odel sentence, are, as just indicated, prohibited by the type theory of \emph{Principia}. \emph{Principia}'s treatment
of Curry's paradox is likewise similar to the type-theory resolution of the liar paradox:
\begin{quote}
  \ldots let $A$ be any sentence, and let $B$ be the sentence:\\
  \indent ``If this sentence is true, then so is $A$.'' \\
  Now we easily see that, if $B$ is true, then so is $A$. That is, $B$ is true. Hence, $A$ is true. We have thus shown that every sentence is true. (\cite{lob1955}: 117)
\end{quote}
To keep matters brief, let's ignore the question of whether truth attaches to sentences or propositions. If we attempt to assign a type to the propositional function predicating truth used in $B$ we see immediately
that any such proposition as $B$ is prohibited by the theory of logical types: if the type of $B$ is $n$, then the type of any propositional function $\phi x$ taking $B$ as an argument must be greater
than $n$ and hence cannot appear as a component of $B$ itself.\footnote{A fuller discussion of the above points could involve an evaluation of G\"odel's  arguments
that examination of proofs of the existence of formally undecidable propositions demonstrate that Whitehead and Russell's solution of the paradoxes is
"too drastic" (\cite{godel1934}:362).}

\section{Conclusion}

The idea of the metamathematical definition of an informal theory, such as the arithmetic of natural numbers, seems initially rather straightforward. When one accepts however that the
associated notion of a formal proof may be diagonalised, using for example G\"odel's approach, doubts may arise as to whether some sort of pathological self-reference might be encountered
in such formal systems. Whilst quite elaborate proofs establishing the consistency of, for example, first-order number theory are available  (\cite{mendelson2015}: Appendix C) it is clear,
in the light of the above, that the notion of a formal proof associated with such systems is not well defined. The metamathematical approach to defining a formal
theory fails to supply reliable devices for the the avoidance of paradox.

If the result is accepted, it follows that metamathematics does not offer a consistent formal theory of the arithmetic of natural numbers.
If the demonstration is reviewed,
it is easy to confirm that the essential requirement for a formal theory to be in scope is that the theory is in scope for G\"odel's proof, hence the
status of a very broad class of formal theories is at issue.
If the conclusion is accepted that a failure of the metamathematical devices for avoiding paradox has been established
the status of metamathematical formal theories that are technically out of scope
is also called into question.

\end{document}